\documentclass[12pt]{article}

\usepackage[all]{xy}

\newtheorem{tm}{Theorem}[subsection]
\newtheorem{lm}[tm]{Lemma}
\newtheorem{pr}[tm]{Proposition}
\newtheorem{rmk}[tm]{Remark}
\newtheorem{cor}[tm]{Corollary}
\newtheorem{ex}[tm]{Example}

\newtheorem{??}[tm]{Question}
\newtheorem{defi}[tm]{Definition}

\newtheorem{ass}[tm]{Assumption}

\font\tenmsb=msbm10
\font\sevenmsb=msbm7
\font\fivemsb=msbm5

\newfam\msbfam
\textfont\msbfam=\tenmsb
\scriptfont\msbfam=\sevenmsb
\scriptscriptfont\msbfam=\fivemsb
\def\Bbb#1{{\fam\msbfam #1}}

\font\teneufm=eufm10
\font\seveneufm=eufm7
\font\fiveeufm=eufm5
\newfam\eufmfam
\textfont\eufmfam=\teneufm
\scriptfont\eufmfam=\seveneufm
\scriptscriptfont\eufmfam=\fiveeufm
\def\frak#1{{\fam\eufmfam\relax#1}}

\def\lorw{\longrightarrow}
\newcommand\n{\noindent}

\newcommand\ci{\cite}

\newcommand\rat{{\Bbb Q}}
\newcommand\comp{{\Bbb C}}
\newcommand\real{{\Bbb R}}
\newcommand\zed{{\Bbb Z}}

\newcommand\pn[1]{{\Bbb P}^{#1}}

\newcommand\blacksquare{{\hspace*{\fill} $\fbox{}$}}

\newcommand{\phix}[2]{ \,^p\!{\cal H}^{#1}(#2)}

\newcommand{\dsdix}[2]{ \bigoplus_{#1} \phix{#1}{#2} [-{#1} ]  }

\newcommand{\coker}{ \hbox{\rm Coker \,} }

\newcommand{\im}{ \hbox{\rm Im \,} }

\newcommand{\ke}{ \hbox{\rm Ker \,} }

\newcommand{\p}{{\Bbb P}}

\newcommand{\ptd}[1]{ \,^{p}\tau_{ \leq {#1} } }

\newcommand{\ptu}[1]{ \,^{p}\tau_{ \geq {#1} } }

\newcommand{\td}[1]{ \tau_{ \leq {#1} } }

\title{Intersection forms, topology of maps and motivic 
decomposition for  resolutions of threefolds}
\author{Mark Andrea A. de Cataldo\thanks{Partially supported by N.S.F.
Grant DMS 0202321}$\;$  
Luca Migliorini\thanks{ Partially supported by
GNSAGA}
}
\date{}

\begin{document}\maketitle

 \tableofcontents

\section{Introduction}
This paper has two aims.

\n
The former is to give an introduction to our earlier work 
\ci{decmightam}
and more generally to some of the main themes of the theory of perverse 
sheaves and
to some of its 
geometric applications. Particular emphasis is put on
the topological properties 
of algebraic maps.  

\n
The latter is to prove 
 a motivic version of the 
decomposition theorem 
for the resolution of a threefold $Y.$ This result allows     
to define  a pure motive whose Betti realization 
is the intersection cohomology of $Y.$

\smallskip
We assume familiarity with Hodge theory and with the formalism 
of derived categories. 
On the other hand, we provide a few explicit computations of perverse 
truncations and intersection cohomology complexes which
we could not 
find in the literature and which  may be helpful 
to understand the machinery.

We discuss in detail the case of surfaces, threefolds and fourfolds. 
In the surface case, 
our ``intersection forms" version of 
the decomposition theorem stems quite naturally from two well-known 
and widely used theorems on surfaces, the  Grauert
contractibility criterion for curves on a surface and the so called 
``Zariski Lemma," 
cf. \ci{BPV}.

\smallskip 
The following assumptions are made throughout the paper

\begin{ass}
\label{s1}
We work with varieties over the complex numbers.
A {\em map} $f:X \to Y$ is a proper morphism of varieties.
We assume that $X$ is smooth. All (co)homology groups are
with rational coefficients.
\end{ass}
These assumptions are placed for ease of exposition only, for
 the main results remain valid when $X$ is singular if one 
replaces the cohomology
of $X$
with its  intersection cohomology, or the constant sheaf
$\rat_{X}$ with the 
intersection cohomology complex of $X.$

\bigskip
It is a pleasure to 
dedicate this work to J. Murre, with admiration and respect.

\section{Intersection forms}
\label{intforms}

\subsection{Surfaces}
\label{surf}
Let $D = \cup D_{k} \subseteq X$ be a finite union of compact 
irreducible curves 
on a smooth complex surface. There is a sequence of
maps
\begin{equation}
\label{blle}
H_{2}(D) \stackrel{r_{*}}\lorw H^{BM}_{2}(X) \stackrel{PD}\simeq
H^{2}(X) \stackrel{r^*}\lorw H^{2}(D).
\end{equation}
The group $H_{2}(D)$ is freely generated by the fundamental classes
$[D_{k}].$

\n
The group $H^{2}(D)\simeq H_{2}(D)^{\vee}$ and, via Mayer-Vietoris,
it is freely generated by the classes associated with points $p_{k} 
\in 
D_{k}.$
The map
$$
H_{2}(D) \stackrel{cl}\lorw H^{2}(X), \qquad cl\,:= \, PD\,\circ \, 
r_{*} 
$$
is called the {\em class map} and it assigns to the fundamental
class $[D_{k}]$ the cohomology class $c_{1}( {\cal O}_{X}(D_{k})  ).$

\n
The restriction map $r,$  or rather
$r\circ PD,$ assigns to a 
Borel Moore $2-$cycle meeting   transversely all the $D_{k},$
the points of intersection with the appropriate multiplicities.

\n
The composition 
$
H_{2}(D) \lorw H^{2}(D),
$
gives rise to the so-called {\rm refined intersection form
on $D \subseteq X:$
\begin{equation}
\label{iota}
\iota: H_{2}(D) \times H_{2}(D) \lorw \rat
\end{equation}
with associated  symmetric intersection  matrix $||D_{h}\cdot 
D_{k}||.$

\n
If $X$  is replaced by  the germ of a neighborhood
of $D,$ then $X$ retracts to $D$ so that
all four spaces appearing in (\ref{blle}) have the same dimension
$b_{2}(D)=$numbers of curves in $D.$

\n
In this case the restriction map 
$r$ is an isomorphism: the Borel Moore classes
of disks transversal to the $D_{k}$ map to the point of intersection.

\n
On the other hand, $cl$ may fail to be injective, e.g.
$(\comp \times \pn{1}, \{0\}\times \pn{1}).$

\medskip
The following are two classical results
results concerning the properties of 
the intersection form $\iota$, dealing respectively with resolutions 
of 
normal surface singularities and one dimensional families of curves.
They are known as the Grauert's Criterion and  the Zariski 
Lemma (cf. \ci{BPV}, p.90).

\begin{tm}
    \label{tmgra}
    Let $f: X \to Y$ be the contraction of a divisor $D$
    to a normal surface singularity.
    Then the refined intersection form $\iota$ on $H_{2}(D)$
    is negative definite.
    In particular, the class map $cl$ is an isomorphism.
\end{tm}

\begin{tm}
    \label{wzar}
    Let $f: X \to Y$ be a surjective proper map of quasi-projective
    smooth varieties, $X$ a surface, $Y$ a curve.
    Let $D = f^{-1}(y)$ be any fiber.
    Then the rank of $cl$ is $b_{2}(D) -1.$
    More precisely, let $F= \sum_k a_k D_{k},$ $a_{k}>0,$ be the
    cycle-theoretic fiber. $F\cdot F=0$ and the induced bilinear form
    $$
    \frac{H_{2}(D)}{ \langle [F] \rangle } \times 
     \frac{H_{2}(D)}{ \langle [F] \rangle }
    \lorw \rat
   $$
is non degenerate and negative definite.
\end{tm}

\begin{rmk}
\label{link}
{\rm Theorem \ref{tmgra} can be interpreted
in terms of the topology of the ``link" ${\cal L}$ of the 
singularity. 
Let $N$ be a small 
contractible neighborhood
of a singular point $y$ and ${\cal L}$ be its boundary.
Choose analytic disks $\Delta_1, \cdots, \Delta_r$ cutting 
transversally 
the divisors $D_1, \cdots, D_r$ at regular points. The classes of 
these disks, 
generate the Borel-Moore homology $H_2^{BM}(f^{-1}(N)) \simeq 
H^2(f^{-1}(N)).$
The statement \ref{tmgra} implies that 
each class $\Delta_{i}$  is 
homologous 
to a  rational linear combination of exceptional curves. 
Equivalently, 
for every index $i$ some multiple of
the $1-$cycle  $\Delta_i\cap {\cal L}$  bounds in the link 
${\cal L}$ of $y.$ 
This is precisely what fails in the aforementioned example
$(\comp \times \pn{1}, \{0\}\times \pn{1}).$
A similar interpretation is possible for the ``Zariski lemma."}

\end{rmk}

In view of the important role played by these theorems 
in the theory of complex surfaces it is natural to ask for 
generalization
to higher dimension.
We next define what is the analogue of the intersection form for a 
general 
map $f:X \to Y$ (cf. \ref{s1})

\subsection{Intersection forms associated to a map}
\label{psif}
General theorems, due to J. Mather, R. Thom and others
(cf. \ci{g-m}) ensure that a projective map
$f:X \to Y$ can be stratified, i.e. there is a decomposition
${\frak Y}= \coprod S_l$ of $Y$ with locally closed nonsingular 
subvarieties
$S_l$, 
the strata, so that $f: f^{-1}(S_l)\to S_l$ is, for any $l$, a 
topologically locally trivial fibration. 
Such stratification allows us, when $X$ is nonsingular, to define
a sequence of intersection forms. 
Let $L$ be the pullback of an ample bundle on $Y.$
The idea is to use sections of $L$ to construct
transverse slices and reduce the strata to points,
and to use a very  ample line bundle $\eta$ on $X$
to fix the ranges: 

Let $ \dim S_l=l$, 
let $s_l$ a generic point of the stratum $S_l$ and  
$Y_s$ a complete intersection of $l$ hyperplane sections of $Y$ 
passing through  $s_l$, transverse to  $S_l$;

As we did for surfaces, we consider the maps:
$$
I_{l,0}: H_{n-l}(f^{-1}(s_l)) \times H_{n-l}(f^{-1}(s_l)) \lorw \rat.
$$ 
obtained intersecting cycles  supported in $f^{-1}(s)$ 
in the smooth $(n-l)-$dimensional ambient variety  $f^{-1}(Y_s):$
$$
H_{n-l}(f^{-1}(s_l))\to H_{n-l}(f^{-1}(Y_s)) \simeq H^{n-l 
}(f^{-1}(Y_s))
\to  H^{n-l}(f^{-1}(s_l)).
$$ 

\medskip
We can define other intersection forms, in different ranges, 
cutting the cycles in $f^{-1}(s_l)$ with generic sections of $\eta.$

The composition: 
$$
H_{n-l -k}(f^{-1}(s))\to H_{n-l -k}(f^{-1}(Y_s))
\simeq H^{n-l +k}(f^{-1}(Y_s)) \to  H^{n-l +k}(f^{-1}(s)).
$$ 
gives maps 
$$
I_{l,k}:H_{n-l-k}(f^{-1}(s_l)) \times H_{n-l+k}(f^{-1}(s_l)) \lorw 
\rat.
$$

Let us denote by 

$$\cap \eta^k:
H_{n-l+k}(f^{-1}(s_l)) \to H_{n-l-k}(f^{-1}(s_l)),
$$ 
the  operation of cutting a cycle in $f^{-1}(s_l)$ with $k$ generic 
sections of $\eta$.

\medskip
\n
Composing this map with $I_{l,k},$ we obtain the intersection forms we 
will consider:
$$
I_{l,k }(\cap \eta^k \cdot, \cdot ): H_{n-l+k}(f^{-1}(s_l)) 
\times H_{n-l+k}(f^{-1}(s_l)) \lorw \rat.
$$ 
\begin{rmk}
\label{inde}
{\rm These intersection forms depend on $\eta$ but not on the 
particular sections 
used to cut the dimension.
They are independent of $L$. 
In fact we could define them using a local slice of the stratum $S_l$ 
and its inverse image, 
without reference to sections of $L.$ }
\end{rmk}

\medskip

\begin{ex}
\label{3fold}
{\rm Let $f:X \to Y$ be a resolution of singularities of a threefold 
$Y,$ 
with a stratification $Y_0 \coprod C \coprod y_0,$
defined so that $f$ is an isomorphism over $Y_0$, 
the fibers are one-dimensional over $C$, and there is a divisor 
$D=\cup D_i$ contracted to the point $y_0$. We have the following 
intersection forms:

\smallskip
-- let $c$ be a general point of $C$ and $s\in H^0(Y, {\cal O}(1))$ 
be a generic section vanishing at $c;$ 
there is the form $H_2(f^{-1}(c) \times  H_2(f^{-1}(c)) \lorw \rat$
which is nothing but the Grauert-type form on the surface 
$f^{-1}(\{s=0\});$

\smallskip
--
similarly, over $y_0$,
there is the form  on $H_4(D)$ given by $\eta \cap [D_i]\cdot [D_j];$ 
it is  a Grauert-type form, 
computed on a hyperplane section 
of $X$ with respect to $\eta;$

\smallskip
--
finally, we have the more interesting  $H_3(D)\times H_3(D) \lorw 
\rat$.
}
\end{ex}

\medskip
One of the dominant themes of this paper is that 
{\em Hodge theory affords non degeneracy results for these forms
and that this non degeneration
has strong cohomological consequences.}
\medskip

To see why Hodge theory is relevant to the study of the intersection 
forms, 
let us sketch a proof of  Theorem \ref{tmgra}, 
in the hypothesis that $X$ and $Y$ are projective. 
The proof we give is certainly not the most natural or economic. 
Its interest lies in the fact that, 
while the original proof seems difficult to generalize to higher 
dimension,
this one can be generalized. 
It is based on the observation that the classes $[D_i]$ of the 
exceptional curves are
``primitive'' with respect to the cup product with the first Chern 
class
of any ample line bundle pulled back from $Y.$ 
Even though such a line bundle is certainly not ample, 
some parts of the ``Hodge package," namely the Hard Lefschetz theorem 
and the 
Hodge Riemann bilinear relations, go through. 
To prove this, we introduce a technique, which we call {\em 
approximation of $L-$primitives,} 
which plays a  decisive role in what follows.

\medskip
\n
{\em Proof of \ref{tmgra} in the case $X$ and $Y$  are projective.}

\n
Let $L$ be the pullback to $X$ of an ample line bundle on $Y$.
Since the map is dominant, $L^2 \neq 0,$ and we get  the 
Hodge-Lefschetz 
type decomposition: 
$$H^2(X,\real) \, = \,  \real \langle c_1(L) \rangle \oplus
\ke \{c_1(L) \wedge : H^2(X) \to H^4(X)\}.$$
Denote the kernel above by $P^{2}.$
This decomposition is orthogonal with respect to the Poincar\'e 
duality 
pairing which,
in turn, 
is non degenerate when restricted to the two summands.
The decomposition holds with rational coefficients. 
However, real coefficients are more convenient in view
of taking limits.

\n
Consider a sequence of Chern classes of ample $\rat -$line bundles 
$L_n$,
converging to the Chern class of $L$, e.g. $L_n=L+\frac{1}{n} \eta,$
$\eta $ ample on $X.$
Define
$P^2_{1/n}=\ke \{c_1(L_n):H^2(X) \to H^4(X)\}.$ 
These are  $(b_2-1)-$dimensional subspaces of $H^2(X)$.
Any limit point of the sequence  $P^2_{1/n}$  in $\p^{b_2}(\real)$ 
gives a codimension one subspace $W \subseteq H^2(X)$,  
contained in $\ke \{c_1(L):H^2(X) \to H^4(X)\}=P^2.$
Since $\dim{W} = b_2-1 = \dim{ P^2},$ we must have $\lim_{n}{ 
P^2_{1/n}}=P^2.$ 

\smallskip
\n
The Hodge Riemann Bilinear Relations hold on $P^2_{1/n}$ by classical 
Hodge theory. 
The duality pairing on the limit
$P^2$ is non degenerate.
It follows that 
the Hodge Riemann Bilinear Relations  hold on $P^2$ as well.

\smallskip
\n
The classes of the exceptional curves $D_i$ are in $P^2,$
since we 
can choose a section of the very ample line bundle on $Y$ not 
passing through the singular point and pull it back to $X.$

\smallskip
\n
The fact that these classes are independent is known
classically. Let us briefly  mention here that if there is only one 
component 
$D_{i}$ then $0 \neq [D_{i}] \in H^{2}(X)$ in the K\"ahler $X.$
In general, one may also argue along the following lines (cf.
\ci{demigsemi}, \ci{deca}, $\S8$): use the Leray spectral sequence
over an affine neighborhood $V$ of the singularity $y$ to show that
$H^{2}(f^{-1}(V)) \to H^{2}(f^{-1}(y))$ is surjective; use the 
basic properties of mixed Hodge structures to deduce
that $H^{2}(X) \to H^{2}(f^{-1}(y))$ is also surjective; conclude
by dualizing and by  Poincar\'e Duality.

\smallskip 
\n

The classes $[D_{i}]$  are real of type $(1,1)$ and for  such 
classes   
 $\alpha\in P^2 \cap H^{1,1}$ the Hodge Riemann
bilinear relations give  
$$
\int_X \alpha \wedge \alpha \,<\, 0 
$$
whence the statement of \ref{tmgra}. 
\blacksquare

\subsection{Resolutions of  isolated singularities in dimension $3$}
\label{lde}
In this section we study the intersection forms 
in the   case  of the resolution of three-dimensional
isolated singularities.  Many  of the features and techniques used
in the general case emerge already in this case.
Besides motivating what follows, we believe 
that the statements and the techniques
used here are of some independent interest.

We prove all the relevant Hodge-theoretic 
results about the intersection forms
associated to the resolution of an isolated singular point on a 
threefold.
This example will be reconsidered in the last section, where we  give 
a motivic version 
of the Hodge theoretic decomposition proved here.

As is  suggested in the proof of Theorem \ref{tmgra} sketched at the 
end of 
the previous section, in order to draw conclusions on the 
behaviour of the intersection forms, we must investigate
the extent to which the Hard Lefschetz theorem and
the Hodge Riemann Bilinear Relations hold when we 
consider the cup product with the Chern class of the 
pullback of an ample bundle by a projective map. 
In order to motivate what follows let us recall an 
inductive proof of the Hard Lefschetz theorem
based on the Hodge Riemann  relations:

\smallskip
\n
{\em 
Hard Lefschetz and Hodge-Riemann relations in dimension 
$(n-1)$ 
and
Weak Lefschetz in dimension $n$ imply Hard Lefschetz in dimension 
$n.$}

\smallskip
\n
Let $X$ be projective nonsingular and $X_H$ be  a generic hyperplane 
section with respect 
to a very ample bundle $\eta$. Consider the map
  $c_1(\eta):H^{n-1}(X)\to H^{n+1}(X).$ 
The Hard Lefschetz theorem 
states it is an isomorphism. By the Weak Lefschetz Theorem
$i^*: H^{n-1}(X)\to H^{n-1}(X_H)$ is injective, and its dual 
$i_*: H^{n-1}(X_H)\to H^{n+1}(X),$ 
with respect to Poincar\'e duality on  $X$ and $X_H$,
is surjective. The  cup product with $c_1(\eta)$ 
is the composition
$i_* \circ i^*$
$$
\xymatrix{
H^{n-1}(X) \ar[rd]^{i^*} \ar[rr]^{c_1(\eta)} &       & H^{n+1}(X) \\
     &     H^{n-1}(X_H) \ar[ur]^{i_*} &
}
$$
and is therefore an isomorphism if and only if the bilinear form
$ \int_{X_{H}}$ 
remains non degenerate when restricted
to the subspace  
$H^{n-1}(X) \subseteq  H^{n-1}(X_H).$ 
This inclusion is a
Hodge substructure.  
The Hodge Riemann relations on $X_{H}$
imply that the Hodge structure $H^{n-1}({X_H})$  is a direct sum of  
Hodge 
structures polarized by the pairing $\int_{X_{H}}.$
It follows that the restriction of the Poincar\'e
form $\int_{X_{H}}$ to $H^{n-1}(X)$ is non degenerate, as wanted. 
The other cases of the Hard Lefschetz Theorem 
(i.e. $c_1(\eta)^k$ for $k \geq 2$) follow immediately from the weak 
Lefschetz theorem and 
the Hard Lefschetz theorem for $X_H$. \blacksquare

\begin{ass}
\label{3foldass}
$Y$ is  projective with an isolated singular point $y$,
$dim Y=3.$
$X$ is a resolution and $f:X \to Y$ is 
an isomorphism when restricted to $f^{-1}(Y-y)$.  
Suppose $D=f^{-1}(y)$ is a divisor and let 
$D_i$ be its irreducible components.
\end{ass}

As usual in this paper, we will denote
by $\eta$ a very ample line bundle on $X,$ 
and by $L$ the pullback to $X$
of a very ample line bundle on $Y.$ 
Of course $L$ is not ample.
We want to investigate whether the Hard Lefschetz theorem and the 
Hodge 
Riemann  relations hold 
if we consider  cup-product with $c_1(L)$ instead of with an ample 
line bundle.

\begin{rmk}
{\rm Since $c_1(L)^3 \neq 0$ we have an isomorphism  
$c_1(L)^3:H^{0}(X) \to H^{6}(X). $}

\end{rmk}

\begin{rmk}
{\rm Clearly the classes $[ D_i ] \in H^2(X)$ are killed by the cup 
product with $c_1(L),$ since
we can pick a generic section of
${\cal O}_Y(1)$ not passing through $y$ and its inverse image in $X$
will not meet the $D_{i}.$
Since $[D_i] \neq 0,$
it follows that  
$c_1(L):H^{2}(X) \to H^{4}(X)$ is not an isomorphism.}
\end{rmk}

We now prove that in fact 
the subspace $\im \{H_4(D) \to H^2(X) \}$ generated by
the classes $[D_{i}]$  is precisely 
$\ke{ c_1(L):H^{2}(X) \to H^{4}(X)}.$

\begin{tm}
\label{ht3fold}
Let $s \in \Gamma(Y,{\cal O}_Y(1))$ be
a generic section and  $X_s=f^{-1}(\{s=0\}) \stackrel{i}{\to}X$.
Then:

a. $i^*:H^1(X) \to H^1(X_s)$ is an isomorphism.

b. $i_*:H^3(X_s) \to H^5(X)$ is an isomorphism.

c. $i^*:H^2(X)/(\im \{H_4(D) \to H^2(X)\}) \to H^2(X_s)$ is injective.

d. $i_*: H^2(X_s) \to \ke \{H^4(X) \to H^4(D)\}$ is surjective.

e. The map $H_3(D) \to H^3(X)$ is injective.
\end{tm}
{\em Proof.}
Set $X_0=X \setminus X_s$ e $Y_o=Y \setminus \{s=0\}$ and let us 
consider 
the Leray spectral sequence
for $f:X_0 \to Y_0.$ Since $Y_0$ is affine, we have $H^k(Y_0)=0$ for 
$k>3$.
$$
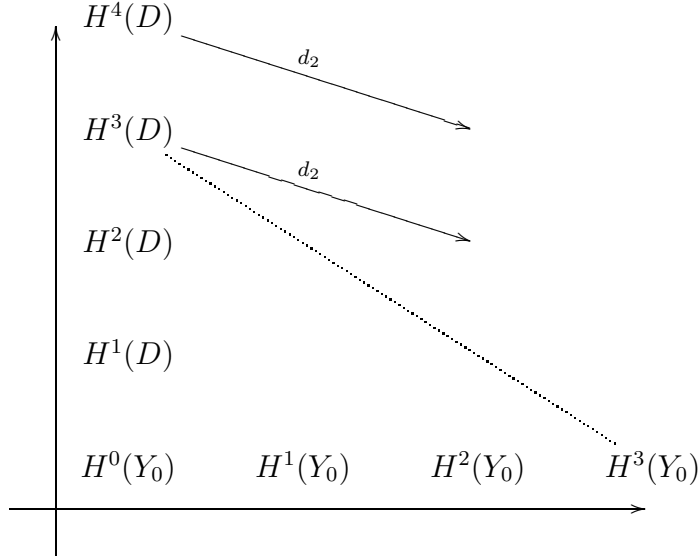
\xymatrix{
  &H^4(D)\ar[rrd]^{d_2}               &         &          &         
&  \\
  &H^3(D)\ar[rrd]^{d_2}\ar@{.}[dddrrr] &         &          &         
&  \\  
  &H^2(D)                              &         &          &         
&  \\
  &H^1(D)                              &         &          &         
&  \\
  & H^0(Y_0)& H^1(Y_0)& H^2(Y_0) & H^3(Y_0)&              \\
\ar@<-4ex>[uuuuu] \ar@<4ex>[rrrr] &                  &          
&         &     }
$$
The sequence degenerates so that we have surjections
$H^3(X_0) \to H^3(D)$ and
$H^4(X_0) \to H^4(D).$ 
But from \ci{ho3}, Proposition 8.2.6,
$H^3(X) \to H^3(D) \to 0$ and
$H^4(X) \to H^4(D) \to 0$ are also surjective.
We  have the long exact sequence
$$
\xymatrix{
H^1_c(X_0) \ar[r]  \ar@{=}[d] & H^1(X) \ar[r] &  H^1(X_s) \ar[rr] 
\ar[dr]& & H^2_c(X_0) \ar[r] \ar@{=}[d]            & H^2(X)  \ar[r]   
& H^2(X_s)   \\  
H^5(X_0)^*=\{0 \}             &               
&                          & 0 \ar[ur] & H_4(X_0) 
\ar@{=}[u]           &        & \\
                              &               &                   
        & & H_4(D) \ar@{=}[u] \ar@{^{(}->}[uur]      &               
& }
        $$
The other statements are obtained applying duality.
\blacksquare

\medskip
Since on $H^2(X_s)$ the bilinear relation of Hodge-Riemann hold,
the argument given at the beginning of this section shows that 
$$
c_1(L)^2:H^1(X)\lorw  H^5(X) \; \hbox{ is an isomorphism}
$$
and
$$
c_1(L)\, :\,  H^2(X)/H_4(D) \, \lorw \,
\ke{ \{H^4(X) \to H^4(D) \}} \; \hbox{ is an isomorphism}.
$$
The Hodge Riemann relations  hold for $P^1:=H^1(X)$, 
$P^2:=\ke c_1(L)^2:H^2(X)/H_4(D) \to H^6(X)  $ 
since, by the weak Lefschetz Theorem, 
they follow from those for $X_s.$ 

\smallskip
The Hodge Riemann relations for  $P^3=\ke \{c_1(L):H^3(X) \to 
H^5(X)\},$ of which $H_{3}(D)$ is a subspace, 
must be considered separately:
the main technique to be used here is the {\em approximation of 
primitives}
introduced in the previous section to prove Theorem \ref{tmgra}.

\begin{tm}
\label{polarizeh3}
The Poincar\'e pairing $\int_X $ is a polarization of $P^3$

\end{tm}
{\em Proof.}
Since 
$c_1(L)^2:H^1(X)\to H^5(X) $ is an isomorphism, there is a 
decomposition,
orthogonal with respect to  the Poincar\`e pairing,
$H^3(X)=P^3 \oplus c_1(L)H^1(X)$ and, in particular,  $\dim 
{P^3}=b_3-b_1, $ just as if
$L$ were ample.
The Poincar\`e pairing remains nondegenerate when restricted to $P^3$.
The classes $c_1(L)+\frac{1}{n}c_1(\eta)$ 
are Chern classes of ample line bundles,
hence
$P^3_{1/n}=\ke \{c_1(L)+\frac{1}{n}c_1(\eta):H^3(X) \to H^5(X) \}$ 
are  $b_3-b_1-$dimensional subspaces of 
$H^3(X)$ .
As in the proof of \ref{tmgra} a limit point of the sequence  
$P^3_{1/n},$ considered as points in the real Grassmannian 
$Gr(b_3-b_1, b_3),$ 
gives a subspace of $H^3(X)$, contained in $\ke \{c_1(L):H^3(X) \to 
H^5(X)\}=P^3$
and, by  equality of  dimensions, $\lim {P^3_{1/n}}=P^3.$ 
The Hodge Riemann relations must then hold 
on the limit $P^3$ as explained in the proof of Theorem
\ref{tmgra}.
\blacksquare

Finally, let us remark 
that the  cup-product with $\eta$ gives an isomorphism $c_1(\eta): 
H_4(D) \to H^4(D)$
via the bilinear form $\int_X c_1(\eta) \wedge [D_i] \wedge [D_j]$ 
which is negative definite. 
As we remarked in \ref{3fold} this form is just 
the intersection form on the exceptional curves of the restriction  
of $f$ to a hyperplane section 
(with respect to  $\eta$) of $X$.

Summarizing:
$$
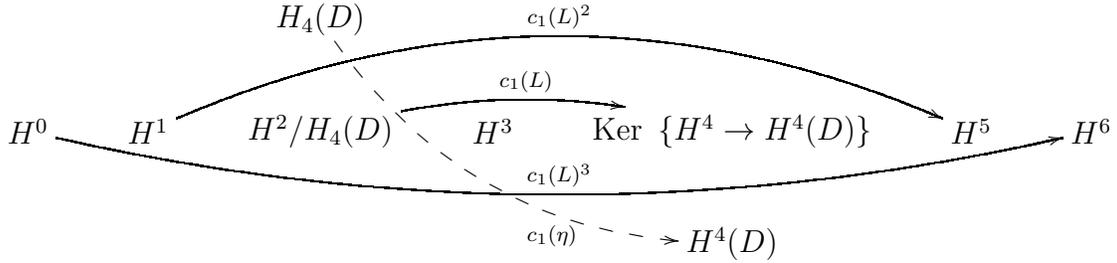
\xymatrix{
         &           &       H_4(D) \ar@/_2pc/ 
@{-->}[ddrr]_>>>>>>>>{c_1(\eta)}    &           
&                                 &               &          \\
H^0  \ar@/_2pc/[rrrrrr]^{c_1(L)^3}  & H^1 
\ar@/^3pc/[rrrr]^{c_1(L)^2}   &  H^2/H_4(D)  \ar@/^1pc/[rr]^{c_1(L)} 
&   H^3  &  \ke \{ H^4 \to H^4(D) \}    &   H^5      &   H^6   \\
         &           &                  &           &     
H^4(D)                      &                &   
}
$$
{\em the groups in the  central row behave, with respect to $L$, 
as the cohomology of a projective nonsingular variety on which $L$ 
is ample.} 

\smallskip
We are now in position to prove the first nontrivial fact on  
intersection forms
which generalizes \ref{tmgra}:
\begin{cor}
\label{tmgra3dim}
$H^3(D)$ has Hodge structure which is pure of weight $3$ and the 
Poincar\'e form is 
a polarization. In particular, the (skew-symmetric) intersection form
$H_3(D) \times H_3(D) \to \rat$ is non degenerate.
\end{cor}
{\em Proof.} This follows because
$H_3(D) \to H^3(X)$, 
is injective and identifies $H^3(D)$ with a  Hodge substructure 
of the polarized Hodge structure  
$ P^3 \subseteq H^3(X)$.
\blacksquare

\medskip
This is a nontrivial 
criterion for a  configuration of (singular) surfaces contained in a 
nonsingular threefold 
to be contractible to a point. See \ci{decmightam}, Corollary 2.1.11
for a generalization to arbitrary dimension.

\n
For example, the purity of the Hodge structure  implies that 
$H^3(D)=\oplus H^3(D_i).$

\smallskip
We will see that the non degeneracy statement of \ref{tmgra3dim}
also plays an important role in the motivic decomposition of $X$
described in section \ref{gmdmt}.

\begin{rmk}
{\rm The same analysis can be carried on with only notational changes 
for an
arbitrary generically finite map 
from a nonsingular threefold $X$, e.g.  assuming that there is 
also some divisor which is blown down to a curve etc. 
In this case the Hodge structure of $X$ can be further decomposed, 
splitting 
of a piece corresponding 
to the contribution to cohomology of this divisor. }
\end{rmk}
 
\begin{rmk}
{\rm The classical argument of Ramanujam \ci{rama}, \ci{e-v},
to derive the Aki\-zu\-ki-Kodaira-Nakano Vanishing Theorem 
from Hodge theory and Weak Lefschetz can be adapted 
to give the following sharp version:
if $L$ is a line bundle on a threefold $X$, 
with $L^3 \neq 0$, 
a multiple of which is globally generated, 
then  
$$
H^{p}(X, \Omega_X^q \otimes L^{-1})=0
$$ 
for $p+q<2$, and for $p+q=2$ but $(p,q) \neq (1,1)$. 
More precisely $H^{1}(X, \Omega_X^1 \otimes L^{-1}) \neq 0$ 
if and only if some divisor is contracted to a point.}
\end{rmk}

\subsection{Resolutions of  isolated singularities in dimension $4$}
\label{lde4}
Let us quickly consider another similar example
in dimension $4$, 

\begin{ass}
\label{4foldass}
$f:X \to Y$,  where $Y$ still has a unique singular point $y$ 
and $X$ is a resolution. As before,  $\eta$ will denote a very ample 
bundle on $X$,
and $L$  the  pull-back of a very ample bundle on  $Y$.
Set  $D=f^{-1}(y).$
\end{ass}
 
An argument completely analogous to the one used in the previous 
example
shows that the sequence of spaces
$H^0(X),$ $H^1(X),$ $ H^2(X)/H_6(D),$ $ H^3(X)/H_5(D),$ $ H^4(X),$ 
$ \ke \{H^5(X) \to H^5(D)\}$, $\ke \{H^6(X) \to H^6(D)\}$, $H^7(X),$ 
$H^8(X) $
satisfies the Hard Lefschetz Theorem with respect to the cup product 
with $L$.
The corresponding primitive spaces $P^1,P^2,P^3,$ are endowed with 
pairing 
satisfying the 
Hodge Riemann bilinear relations.
The new fact that we have to face shows up when studying the Hodge 
Riemann 
bilinear relations
on $H^4(X)$. 
The ``approximation of primitives'' technique here must be modified, 
since
the dimension of $P^4= \ke \{c_1(L):H^4(X) \to H^6(X) \}$ is greater 
than
$b_4-b_2. $ Hence, if we introduce the primitive spaces 
$P^4_{1/n}=\ke \{c_1(L)+\frac{1}{n}c_1(\eta):H^4(X) \to H^6(X)\}$ 
with respect to the ample classes $c_1(L)+\frac{1}{n}c_1(\eta),$ 
their limit is a proper subspace, of dimension $b_4-b_2$, of $P^4$. 
We can determine the exact dimension of $P^4$:

\begin{lm}
$\dim{ \ke \{c_1(L):H^4(X) \to H^6(X)\}}\,=\, b_4  -b_2 + 
\dim{H_6(D)}.$
\end{lm}
{\em Proof.}
Since $c_1(L)^2:H^2(X)/H_6(D) \stackrel{c_1(L)}{\to} 
H^4(X) \stackrel{c_1(L)}{\to}\ke \{H^6(X) \to H^6(D)\}$
is an isomorphism, we have an orthogonal decomposition
$$
H^4(X)= P^4 \oplus \im \{c_1(L): H^2(X) \to H^4(X)\}.
$$
The statement follows from:
$ \ke \{c_1(L): H^2(X) \to H^4(X) \}=H_6(D).$
\blacksquare

\medskip
The ``excess'' dimension of $P_4$ is thus 
$\dim{H_6(D)}$. On the other hand  $P^4$
contains an obvious subspace of this dimension, namely 
$c_1(\eta)H_6(D)$, the subspace generated by the classes obtained 
intersecting
the irreducible components of the exceptional divisor with a generic  
hyperplane section.

\begin{rmk}
{\rm The intersection form $\int_X c_1(\eta)^2 \wedge [D_i] \wedge 
[D_j]$ 
is negative definite, 
as it is just the intersection form on the exceptional curves of a 
double
hyperplane section
of $X$.}
\end{rmk}

This last remark implies the following  orthogonal decomposition
$$
H^4(X)= c_1(\eta)H_6(D) \oplus  (c_1(\eta)H_6(D))^{\perp}\cap P^4 
\oplus \im \{  c_1(L): H^2(X) \to H^4(X) \}.
$$
and  $(c_1(\eta)H_6(D))^{\perp}\cap P^4$ has dimension $b_4-b_2.$
This subspace turns out to be the subspace of ``approximable  
$L-$primitives'' 
we are looking for, as shown in the following

\begin{tm}
$$
\lim_{n \to \infty}{
\ke{ \{ c_1(L)+\frac{1}{n}c_1(\eta) :H^4(X) \to H^6(X)\}}
} \, 
=\, (c_1(\eta)\,H_6(D))^{\perp}\cap P^4.
$$
\end{tm}
{\em Proof.} The two subspaces have the same dimension, so it is 
enough 
to prove that 
$$ \ke{ \{
c_1(L)+\frac{1}{n}c_1(\eta): \{ H^4(X) \to H^6(X) \}
\}
}
\subseteq (c_1(\eta)H_6(D))^{\perp}.
$$
If $(c_1(L)+\frac{1}{n}c_1(\eta)(\alpha))=0,$ then, using
$c_{1}(L) [D_{i}]=0:$

$$
\int_X c_1(\eta)\wedge [D_i] \wedge \alpha \, = \, 
-n \int_X  c_1(L) \wedge [D_i] \wedge \alpha=0.
$$
\blacksquare

\begin{cor}
The Poincar\'e pairing is a polarization of the weight $4$ pure  
Hodge structure 
 $(c_1(\eta)H_6(D))^{\perp}\cap P^4$.
\end{cor}

\medskip
Let us spell out the consequences of this analysis for the 
intersection form
$ H_4(D) \times H_4(D) \to \rat.$ 
First notice that the same argument used in the proof of  
\ref{ht3fold}.e shows that the map 
$H_4(D) \to H^4(X)$ is injective. It follows that
$H_4(D) $ has a pure Hodge structure
 which is the direct sum  of two substructures 
polarized (with opposite signs) by the  Poincar\'e pairing.

The next result shows that in fact $H_4(D) $ is the direct sum  
of two substructures, polarized (with opposite signs) by the  
Poincar\'e pairing.
This gives a clear indication of what happens in 
general:

\begin{cor}
\label{tmgra4fold}
The intersection form $ H_4(D) \times H_4(D) \to \rat$ is 
non degenerate.
There is a direct sum decomposition:
$$
H_4(D)= c_1(\eta)H_6(D) \oplus (c_1(\eta)H_6(D))^{\perp}
$$
orthogonal with respect to the intersection form, which is negative 
definite on
the first summand and positive on the second.  
\end{cor}

\section{Intersection forms and Decomposition in the Derived Category}
\label{intformdecomp}
We now show how the results we quoted at the beginning of the first 
section
can be translated in statements about 
the decomposition in the derived category of sheaves 
of the direct image of the constant sheaf.
We will freely use the language of derived categories.
In particular we will  use the notion of a constructible sheaf and 
the functors 
$Rf_*,Rf_!,f^*,f^!.$

In section \ref{subsm}
we  briefly review the classical $E_2-$degeneration criterion
of Deligne \ci{dess}, \ci{shockwave} in order to motivate 
the construction
of the perverse cohomology complexes.
These complexes are
 a natural generalization 
of the higher direct image local systems for a smooth map.
The construction of perverse cohomology is carried out in section 
\ref{macc}.

\medskip
We denote by $S(Y)$ the abelian category of sheaves of $\rat -$vector 
spaces
on $Y$, and by $D^b(Y)$ the corresponding derived category of 
bounded complexes.

\medskip
We shall make use of the following splitting criterion
in the derived category. We state it in the form we need it in this 
paper. For a more general statement and a proof the reader is 
referred to \ci{demigsemi} and \ci{decmightam}.

Let $(U,y)$  be a germ of an isolated $n-$dimensional singularity 
with the 
obvious stratification  $= V \coprod y,$ $j: V \to U \leftarrow y : 
i$ be the 
obvious maps, $P$ be a self-dual
complex on $U$ with $P_{|V}= {\cal L}[n],$ ${\cal L}$
a local system on $V,$ and $P \simeq \tau_{\leq 0}P. $
We wish to compare
$P,$ $IC_{U}({\cal L}):= \td{-1} Rj_{*} {\cal L}[n]$
and the stalk ${\cal H}^{0}(P)_{y}.$

\begin{lm}
\label{splitp}
The following are equivalent:

\smallskip
\n
1) there is a canonical isomorphism  in the derived category 
$$
P \, \simeq \,  IC_{U}({\cal L}) \oplus
{\cal H}^{0}(P)_{y}[0];
$$

\n
2) the natural map ${\cal H}^{0}(P) \lorw {\cal H}^{0}(Rj_{*} 
j^{*}P)=R^{n}j_{*}{\cal L}$ 
is zero.
\end{lm}

\subsection{Resolution of surface singularities }
\label{fsts}
For a normal surface  $Y$, let $j:Y_{reg} \to Y$ be the open 
embedding 
of its regular points. The intersection cohomology complex, 
which we will consider in much more detail in the next section, is 
$IC_Y= \tau_{\leq -1} Rj_{*}\rat_{Y_{reg} } [2] .$
The following, which we will prove as a consequence of \ref{tmgra}, 
is the first case of the Decomposition theorem which 
needs to be stated in the derived category and not  just 
in the category of sheaves.

\begin{tm}
\label{tmrs}
Let $f:X \to Y$  be a proper birational map of
quasi-projective surfaces, $X$ smooth, $Y$ normal.
There is a  canonical isomorphisms
$$
Rf_{*} \rat_{X}[2] \stackrel{\simeq}\lorw IC_{Y} \oplus 
R^{2}f_*\rat_X [0].
$$
\end{tm}
{\em Proof.}
We work locally on $Y.$ 
Let $(Y,p)$ be the germ of an analytic normal surface singularity,
$f:(X,D) \to (Y,p)$ be a resolution. The fiber $D=f^{-1}(p)$ is a 
connected union of 
finitely many irreducible compact curves $D_k.$
Note that $j^*Rf_*\rat \simeq \rat_{Y\setminus p}.$
Consider the following diagram
$$
  \xymatrix{
Rf_{*}\rat_{X}[2]   \ar[rd]_{l_{0}} \ar@/_1pc/[rdd]_{l_{-1}} 
\ar@/_2pc/[rddd]_{l_{-2}}  
\ar[r]^{r} &
Rj_{*}j^{*} Rf_{*}\rat_{X}[2] & =  Rj_{*} \rat_{Y 
\setminus p}[2] &   \\   
&  \td{0} Rj_{*}j^{*} Rf_{*}\rat_{X}[2] \ar@[=][u] 
\ar@[=][u] &  &   \\
  & \td{-1} Rj_{*}j^{*} Rf_{*}\rat_{X}[2]   \ar[u]  & 
  =:  IC_{Y} &\\
 &  \td{-2} Rj_{*}j^{*} Rf_{*}\rat_{X}[2] \ar[u] & =   
 \rat_{Y}[2] &
}
$$
We are   looking for the obstructions to the existence
of the lifts $l_{0},$ $l_{-1}$ and $l_{-2}.$

\n
Since $R^k f_*\rat_X =0,$ $k\geq 3,$
we have that  $\td{0} Rf_*\rat_X[2] \simeq Rf_* \rat_X[2].$
In particular, $l_0$ exists and it is unique.

\n
From the exact triangle
$$
 \to \tau_{\leq -1}Rj_*j^*Rf_* \rat_X[2]
\to  \tau_{\leq 0}Rj_*j^*Rf_* \rat_X[2]
\to {\cal H}^0 (Rj_*j^*Rf_* \rat_X[2])\simeq R^2j_* \rat_{Y \setminus 
p} 
\stackrel{+1}{\to}
$$

$l_{-1}$ exists iff
the natural map
$$
\rho: R^2f_* \rat_X \to R^2j_* \rat_{Y \setminus p}
$$
is trivial. Using the isomorphisms $(R^2f_* \rat_X)_{p} \simeq
H^{2}(X)$ and $(R^2j_* \rat_{Y \setminus p})_{p} \simeq
H^{2}(X \setminus D),$ the map $\rho$
can be identified with the restriction map $\rho$ appearing
in the long exact sequence of the
pair $(X, D):$
$$
\ldots \lorw  H_2(D) \stackrel{cl}\lorw H^2(X)\simeq H^2(D) 
\stackrel{\rho}\lorw
H^2( X \setminus D ) \lorw \ldots 
$$
where  have identified $H_{2}(D) \simeq H^{2}(X, X \setminus D)$
via Lefschetz Duality.

\n
By Theorem \ref{tmgra}, $cl$ is an isomorphism
and $\rho$ is trivial.

\n
This lift $l_{-1}$ is unique and splits by 
Lemma \ref{splitp}.
\blacksquare

\begin{rmk}
\label{nolift}
{\rm 
It can be shown easily that a lift 
$Rf_{*}\rat_{X}[2] \lorw \td{-2}\, Rj_{*}\rat_{U}[2] = \rat_{Y}[2]$
exists iff $H^{1}(f^{-1}(p), \rat) =\{0\},$ i.e. iff $IC_{Y}\simeq
\rat_{Y}[2],$ iff $(Y,p)$ is a rational homology manifold.

\n
It follows that, in general, the natural map $\rat_{Y } \to 
Rf_{*}\rat_{X}$ 
does {\em not}  split and $Rf_{*}\rat_{X}$  does {\em not}
 decompose as a direct sum of its shifted cohomology sheaves
 as in (\ref{e1}).
 }
\end{rmk}

\begin{ex}
    \label{rag2}
    {\rm
Let $f: X = \comp \times \pn{1} \to Y$ be the real algebraic map
contracting  precisely $D:= \{0 \} \times \pn{1}$ to a point $p \in 
Y.$
One has 
a {\em non} split  exact sequence in the category 
$P(Y)$ of perverse sheaves 
on 
$Y:$
$$
 0  \lorw IC_{Y} \lorw Rf_{*}\rat_{X}[2] \lorw H^{2}(\pn{1})_{p}[0] 
 \lorw 0.
$$
}
\end{ex}

\bigskip
It is remarkable that while the lift $l_{-2}$ does not exist in 
general,
the lift $l_{-1}$ always exists.
While looking for a nontrivial map $Rf_{*} \rat_{X}[2] \to 
\rat_{Y}[2],$ 
one  ends up  finding another    {\em more interesting} map to 
$IC_{Y}.$

\bigskip
Recall that the dualizing sheaf $\omega_{X} \simeq \rat_{X}[2n].$
Dualizing the canonical isomorphism of Theorem \ref{tmrs}
and, keeping in mind that $IC_{Y}$ and $H_{2}(D)_{p}[0]$
are simple objects in $P(Y)$ (cf. section \ref{simple}), we get
\begin{cor}
\label{ctmrs} There are canonical isomorphisms
$$
H_{2}(D)_{p}[0] \oplus IC_{Y}^{*} \stackrel{\simeq}\lorw
Rf_{*} \omega_{X}[-2] \stackrel{PD}\simeq
Rf_{*} \rat_{X}[2] \stackrel{\simeq}\lorw IC_{Y} \oplus 
H^{2}(D)_{p}[0],
$$
such that the composition is a direct sum map
and induces the intersection form $\iota$ on $H_{2}(D)$ and the 
Poincar\'e-Verdier pairing on the self-dual $IC_{Y}.$

\n
In particular, if $X$ is compact, then the induced splitting 
injection 
$$
I\!H^{\bullet}(Y) \subseteq H^{\bullet}(X) 
$$
exhibits the lhs as the pure Hodge substructure of the rhs
orthogonal to the space $cl (H_{2}(D)) \subseteq H^{2}(X)$
with respect to the Poincar\'e pairing on $X.$
\end{cor}

\subsection{Fibrations over  curves}
\label{fstc}
Let $f: X \to Y$ be a map of from a smooth surface onto a smooth 
curve.
Denote by $\hat{f}: \hat{X} \to \hat{Y}$ 
the smooth part of the map $f,$ 
by $j: \hat{Y} \to Y$ the open immersion,
by  $T^i:= R^i\hat{f}_*\rat_{\hat{X}} = {R^if_*\rat_X}_{|\hat{Y}}.$
For ease of exposition we assume that $f$ has connected fibers.

\smallskip
\n
Fix an ample line bundle $\eta$ on $X.$ The isomorphism
stated in the next proposition will depend on $\eta.$

\begin{pr}
\label{rhl1}
There is an isomorphism
$$
Rf_* \rat_X[2] \simeq j_*T^0[2] \oplus P
\oplus j_*T^2[0],
$$
with $P$ a suitable self-dual (with respect to the Verdier duality
functor) object of $D^b(Y).$
\end{pr}
{\em Proof.}
We  work around one critical value $p \in Y$
and replace $Y$ by a small disk centered at $p,$ $X$ by the preimage 
of this disk, etc.

\n
Since the fibers are connected, $\rat_Y \simeq f_*\rat_X \simeq 
j_*T^0 \simeq
j_*T^2.$ 

\n
Since $\eta$ is $f-$ample, $\eta: T^0 \simeq T^2,$
which in this case implies that $\eta: j_*T^0 \simeq j_*T^2.$

\n
There are the natural truncation 
maps $f_*\rat_X \to Rf_*\rat_X \to R^2f_*\rat_X[-2].$

\n
There is the natural adjunction map $R^2f_*\rat_X \to j_*T^2.$
It is splitting/surjective in view of the presence of $\eta.$

\n
Putting together, 
there is a sequence of maps
\begin{equation}
\label{som}
j_*T^0[2] \stackrel{c}\to Rf_*\rat_X[2] \stackrel{\eta}\to Rf_* 
\rat_X[4]
\stackrel{\pi}\to j_*T^2[2] \stackrel{\sigma}\to
j_*T^0[2] \stackrel{c}\to \ldots
\end{equation}
where the composition $\eta \pi c$ is the isomorphism
mentioned above
$\eta: j_*T^0 \simeq j_*T^2$ and $\sigma := ( \eta \pi c  )^{-1}.$

\n
The reader can verify that the composition
$$
j_*T^0[2] \oplus j_*T^2[0] \stackrel{
\gamma:=( c + \eta c \sigma)[-2]     }\lorw 
Rf_* \rat_X[2] 
\stackrel{ \sigma\pi\eta \oplus \pi [-2]}\lorw j_*T^0[2] \oplus 
j_*T^2[0]
$$
is the identity, i.e. $\gamma$ splits.

\n
Let $P:= \mbox{Cone}(\gamma).$ There is a direct sum decomposition
\begin{equation}
\label{chsp}
Rf_* \rat_X [2] \simeq j_*T^0[2] \oplus P \oplus j_*T^2[0].
\end{equation}
The self-duality of $P$ follows from the self-duality of 
$Rf_* \rat_X [2]$
and of $j_*T^0[2] \oplus j_*T^2[0].$
\blacksquare

\medskip
The object $P$ 
introduced in the previous proposition has a simple structure:
\begin{pr}
\label{semisimp1}
Assumptions as in  \ref{rhl1}. 
The object $P$ splits in $D^b(Y)$ as
$$P=V \oplus j_*T^1[1],$$
where $V = \ke{R^2f_*\rat_X \to j_*T^2}$
 is a  skyscraper sheaf supported at
$Y \setminus \hat{Y}.$ 
This decomposition is canonical and compatible with Verdier duality.
\end{pr}
{\em Proof.}
By inspecting cohomology sheaves we  see that
${\cal H}^i(P) =0$ for $i \neq 0,1,$ that 
${\cal H}^{-1}(P) = R^1f_*\rat_X$ and that
${\cal H}^0(P) =V.$

\n
In view of Lemma \ref{splitp}, we need to show that
\begin{equation}
    \label{rpr}
r': {\cal H}^0 (P) \to R^1j_*T^1
\end{equation}
is the zero map.

\n
We now show that this is equivalent to the Zariski Lemma.

\n
By applying adjunction to (\ref{chsp}), we obtain the  a commutative 
diagram
$$
  \xymatrix{
Rf_{*}\rat_{X}[2]   \ar[d]^{\simeq}   
\ar[r] &
Rj_{*}j^{*} Rf_{*}\rat_{X}[2] \ar[d]^{\simeq}     \\   
j_{*}T^{0}[2] \oplus P \oplus j_{*}T^{2}[0] \ar[r] &
Rj_{*}T^{0}[2] \oplus Rj_{*}P \oplus Rj_{*}T^{2}[0]
}
$$
The associated map of spectral sequences
${\Bbb H}^p(Y, {\cal H}^q (-) ) \Longrightarrow
{\Bbb H}^{p+q}(Y, -)$
gives a commutative diagram 

$$
  \xymatrix{
H_{2}(D) \ar[r]^{cl} & H^{2}(X)   \ar[d]^{\simeq}   
\ar[r]^{r} &
H^{2}(X \setminus D) \ar[d]^{\simeq}     \\   
& {\cal H}^{0}( P)  \oplus (j_{*}T^{2})_{p} \ar[r]^{(r',id)} & \;
{\cal H}^{1}( Rj_{*}T^{1})  \oplus (j_{*}T^{2})_{p}. 
}
$$
It follows that, using the identifications above, 
$\mbox{Im}(cl)=\ke{r} = \ke{r'}
\subseteq {\cal H}^0(P).$

\n
In particular, $r'=0,$ iff $\dim{ \ke{r} } = \dim{ {\cal H}^0(P) }.$
Note
 that $\dim{ {\cal H}^0(P) }=b_2(f^{-1}(p)) -1.$

\n
It follows that $r'=0$ iff $\dim{ \mbox{Im} (cl) } = 
b_2(f^{-1}(p)) -1.$ The latter is implied by 
Theorem \ref{wzar}.

\n
We conclude by  Lemma \ref{splitp}.
\blacksquare

\medskip
Finally, we have:

\begin{tm}
\label{stc} 
There is an isomorphism
$$
Rf_* \rat_X[2] \simeq j_*T^0[2] \oplus j_*T^1[1]  \oplus V [0]
\oplus j_*T^2[0].
$$
\end{tm}

\begin{rmk}
{\rm From \ref{stc} follows that $R^1f_* \rat_X \simeq 
j_*R^1\hat{f}_*\rat_{\hat{X}}.$
Note that this implies the Local Invariant Cycle Theorem.

\n
Since $R^2f_* \rat_X=j_*R^2\hat{f}_*\rat_{\hat{X}}\oplus V,$ we have 
the coarser decomposition 
$$
Rf_* \rat_X[2] \simeq R^0f_*\rat_{X}[2] \oplus R^1f_*\rat_{X}[1]
\oplus R^2 f_*\rat_{X}[0].
$$
In particular, the Leray spectral sequence degenerates at $E_2.$
It is easy to see that the Leray filtration on the cohomology of $X$ 
is by Hodge substructures:
$$
L^2=f^*H(Y), \qquad L^1= \ke{ \{f_*:H(X) \to H(Y)\}}.
$$ }
\end{rmk}

\subsection{Smooth maps}
\label{subsm}
Even in the case of a smooth fibration $f: X\to Y$
of a surface over a curve, the study of the complex
$Rf_{*}\rat_{X}$ is nontrivial, without any projectivity assumptions.

\begin{ex}
\label{hopf}
{\rm
Let $X$ be a Hopf surface.
There is a  natural holomorphic smooth
fibration $ f:X\to \pn{1}$ with fibers elliptic curves.
Since $b_{1}(X)=1,$ one sees easily that the 
Leray Spectral Sequence for $f$ is not $E_{2}-$degenerate.
In particular, $Rf_{*}\rat_{X} $ is {\em not}
isomorphic to $\oplus_{i}{ R^{i}f_{*}\rat_{X}  
[-i]}.$
}
\end{ex}

\medskip
Let us briefly list {\em some} 
of the important properties of  a smooth projective  map 
$f:X \to Y$ of 
smooth varieties.

\medskip
The sheaves $R^{i}f_{*}\rat_{X}$ are locally constant over $Y,$
i.e. they are local systems.
In fact, $f$ is differentiably locally trivial over $Y$
in view of   Ehresmann Lemma.

\medskip
The model for a general decomposition theorem for $Rf_{*}\rat_{X}$
is the following
\begin{tm}
Let      $f:X^{n} \to Y^{m}$ be a smooth projective map
of smooth quasi-projective varieties
of the indicated dimensions and  $\eta$  be an $f-$ample
line bundle on $X.$
    \label{del}
\begin{equation}
    \label{e1}
    Rf_{*}\rat_{X} \, \simeq_{D(Y)}  \, \oplus{R^{i}f_{*}\rat_{X} 
[-i]}.
\end{equation}

\begin{equation}
    \label{e2}
\eta^{i}: R^{n-m-i}f_{*}\rat_{X} \simeq 
R^{n-m+i}f_{*}\rat_{X}, \quad  \forall \,i \geq 0.
\end{equation}

\begin{equation}
    \label{e3}
 \mbox{   The local systems $R^{j}f_{*}\rat_{X}$ on $Y$
are semisimple.}
\end{equation}
\end{tm}

\medskip
The  first Chern class of the
line bundle $\eta \in H^{2}(X, \rat) =
Hom_{D(X)}(\rat_{X}, \rat[2]),$ defines maps 
$$
\eta: \rat_{X} \lorw 
\rat_{X}[2],   \quad \eta: Rf_{*}\rat_{X} \lorw Rf_{*}\rat_{X}[2],
\quad
\eta^{r}:  Rf_{*}\rat_{X} \lorw Rf_{*}\rat_{X}[2r]
$$
and finally
$$
\eta^{r}: R^{i}f_{*}\rat_{X} \lorw R^{i+2r}f_{*}\rat_{X}.
$$
Theorem \ref{del}.\ref{e2} is then just a re-formulation
of the Hard Lefschetz Theorem for the fibers of $f$
and can be named the Relative Hard Lefschetz Theorem
for smooth maps.

\medskip

We remind the reader that a functor 
${\cal T}: D(Y) \to A$, $A$ an abelian category,
is said  to be cohomological (cf. \ci{verd} II, 1.1.5.) if, setting 
${\cal T}^i(K)= {\cal T}^0(K[i]),$  to a distinguished triangle
$$
K \to L \to M \stackrel{+1}{\to}
$$
corresponds a long exact sequence in $A$

$$
\to {\cal T}^i(K) \to {\cal T}^i(L) \to {\cal T}^i(M) \to {\cal 
T}^{i+1}(K) \to... 
$$
The cohomology sheaf functor ${\cal H}^{0}: D(Y) \to S(Y)$
is cohomological. 
Noting that ${\cal H}^{i}(Rf_{*}) = R^{i}f_{*},$
Theorem \ref{del}.\ref{e1} can be re-phrased by saying
that $Rf_{*}\rat_{X}$ is decomposable with respect to
the functor ${\cal H}^{0}.$ 

\n
It is important to note that (\ref{e2}) implies (\ref{e1}) by
the Deligne-Lefschetz Criterion.

\n
Theorem \ref{del}.\ref{e3}
states that every local subsystem  ${\cal L} \subseteq
R^{j}f_{*}\rat_{X}$ admits a complement, i.e. a local system ${\cal 
L}'$
such that ${\cal L} \oplus {\cal L}' = R^{j}f_{*}\rat_{X}.$

\medskip
Let us note {\em some} of the important consequences of Theorem
\ref{del}.

\n
The ${\cal H}^{0}-$decomposability (\ref{e1}) of $Rf_{*}\rat_{X}$
implies immediately the $E_{2}-$degeneration of the Leray spectral 
sequence,
i.e. of the spectral sequence associated with the
cohomological functor ${\cal H}^{0}:$
\begin{equation}
    \label{e4}
\mbox{ $H^{p}(Y, R^{q}f_{*}\rat_{X}) \Longrightarrow  H^{p+q}(X, 
\rat)$ is $E_{2}-$degenerate}.
\end{equation}
This degeneration implies the surjection
\begin{equation}
    \label{e5a}
H^{k}(X, \rat) \lorw H^{0}(Y, R^{k}f_{*}\rat_{X}) = H^{k}(X_{y}, 
\rat)^{\pi_{1}(Y, y)},
\end{equation}
i.e. the  so-called Global Invariant Cycle Theorem.

\n
The Theory of MHS allows to show, using a smooth
compactification of $X,$ that in fact the
monodromy invariants are a Hodge substructure of $H^{k}(X_{y},\rat),$
which as a PHS is independent of $y \in Y$ (Theorem of the Fixed 
Part).
 In fact, (\ref{e3}) is a consequence of this fact.

\bigskip
In general, if $f$ is not smooth, Theorem \ref{del}
fails completely.

\n
The Relative Hard Lefschetz Theorem (\ref{e2}) fails 
due to the presence of singular fibers, i.e.
fibers along which
the differential of $f$ drops rank.

\n
The sheaves $R^{j}f_{*}\rat_{X}$ are no longer locally constant.
Moreover, they are not semisimple in the category of constructible 
sheaves: e.g. $j_{!}\rat_{\comp^{*}} \to \rat_{\comp}.$

\n
The following examples shows that
the ${\cal H}^{0}-$decomposability (\ref{e1}) fails in general and so 
does 
the $E_{2}-$degeneration of the Leray Spectral Sequence 
(\ref{e4}).

\begin{ex}
    \label{exnondeg}
    {\rm
Let $X$ be the blowing up of $\comp \pn{2},$
along ten
points lying on an irreducible cubic $C'$ and $C$ be the strict
transform of $C'$ on $X.$ Since $C^{2}=-1$ the curve contracts
to a point under a birational map $f: X\to Y.$
We leave to the reader the task to verify that
1) the Leray Spectral Sequences  for $H^{2}(X,\rat)$
and for $I\!H^{2}(Y, \rat)$ 
are not $E_{2}-$degenerate and that, though
the Leray spectral sequence always degenerates over suitably small 
Euclidean neighborhoods on $Y,$
2) the complex $IC_{Y}$ does not split
as a direct sum of its shifted cohomology sheaves.
}
\end{ex}

\medskip
The following more general class of examples shows that  the failure
of the $E_{2}-$degeneration is very frequent.
\begin{ex}
    \label{exgen}
{\rm
Let $f: X \to Y$ be a  projective resolution of the singularities
of a projective  and normal variety $Y$ such that
there is at least one index 
$i$ such that the  natural MHS 
on  $H^{i}(Y, \rat)$ is not pure (e.g. $i=2$ in\ref{exnondeg}). 
Then the Leray Spectral Sequence for $f$ is not
$E_{2}-$degenerate. 
If it were, then the edge sequence would
give an injection of MHS  $f^{*}: H^{j}(Y, \rat)
\to H^{j}(X, \rat),$  forcing
the MHS of such a $Y$ to be pure.
}
\end{ex}

\medskip
However, not everything is lost.

\section{Perverse sheaves and the Decomposition Theorem}
\label{macc}
One of the main 
ideas  leading to the theory of perverse 
sheaves  is that Theorem \ref{del}, which holds for smooth maps,
can be made to hold for arbitrary
proper algebraic maps  provided that  it is re-formulated using
the perverse cohomology functor  $^{p}\!{\cal H}^{0}$ in place of the
cohomology sheaf functor ${\cal H}^{0}.$
Just as this latter is $\tau_{\leq 0}\tau_{\geq 0},$ with 
$\tau$ the standard truncation functors of a complex,
 the perverse cohomology functor will be expressed
as $^{p}\!{\cal H}^{0}= \, ^{p}\!\tau_{\leq 0} \, \!^{p}\!\tau_{\geq 
0},$
where $^{p}\!\tau$ is the so called perverse truncation functor. 
Roughly speaking, 
the perverse truncation functor (with respect to middle perversity, 
which is the only case we will consider)
is defined by gluing standard truncations on the strata,
shifted by a term which depends on the dimension of the stratum.
The choice of the shifting is dictated by the behavior of the 
standard 
truncation with respect to duality,
as we suggest in \ref{troncodiaz}.
In this context and keeping this in mind,  
perverse truncation becomes quite natural.
We believe it can be 
useful to give a few details of its construction and an example of 
computation, 
related to the examples given in  section
\ref{psif}. In analogy with the cohomology sheaf functor ${\cal 
H}^0,$ 
the perverse cohomology functor $^{p}\!{\cal H}^{0}$ will be a 
cohomological 
functor which 
takes values in an abelian 
subcategory of ${D}^{b}(Y),$
whose object are the so called perverse sheaves. For a general proper 
map these 
objects  play the role 
played by local system for smooth maps.

\subsection{Truncation and Perverse sheaves}
\label{sps}
Let ${D}^{b}(Y)$ be the  bounded derived category of the category
$S(Y)$
of  sheaves of rational vector spaces on $Y.$
We are interested in the full subcategory $D(Y)$ of those complexes
whose cohomology sheaves are constructible. This means that,
given an object $F$ of $D(Y),$ there is an algebraic Whitney 
stratification $Y=\coprod{S_{l}},$ depending on $F,$  such that
${\cal H}^{j}(F)_{|S_{l}}$ is a finite rank local system.
By the Thom Isotopy Lemmata, $Rf_{*}\rat_{X}$,
and in fact any other complex appearing in this paper, is
an object of $D(Y).$ One is interested in direct sum decompositions
of this complex, in the   geometric meaning of the summands
and in the consequences, both theoretical and practical,
of such splittings.

We now define the $t$-structure on $D(Y)$ associated with the middle 
perversity. Instead of insisting on its axiomatic characterization 
(cf. 
\ci{bbd}),
we give the explicit construction of 
the {\em perverse truncations}
     $\ptd{m}: D(Y) \lorw D(Y)$,
and     $\ptu{m}: D(Y) \lorw D(Y)$.
These come with natural morphisms
 $\ptd{m}F \lorw F$ and
 $F \lorw  \ptu{m}F.$

\n
We start with the following:
\begin{lm}
\label{troncodiaz}
Let $Z$ be nonsingular of complex dimension $r$, and $F \in D(Z)$ 
with locally constant cohomology sheaves. Then there are
natural isomorphisms:
$$
\tau_{\leq k} {\cal D}F \simeq  {\cal D}\tau_{\geq -k-2r}F
\qquad 
\tau_{\geq k} {\cal D}F \simeq  {\cal D}\tau_{\leq -k-2r}F
$$
\end{lm}
{\em Proof.}
Since the dualizing complex is in this case isomorphic to
$\rat_Z[2r],$ it is enough to prove that there are natural 
isomorphisms
$$
\tau_{\leq k} Rhom(F,\rat_Z) \simeq  Rhom(\tau_{\geq -k}F, \rat_Z)
\qquad 
\tau_{\geq k} Rhom(F,\rat_Z) \simeq  Rhom(\tau_{\leq -k}F, \rat_Z).
$$
We prove the first statement. The proof of the 
second is analogous.
Applying $Rhom$ and $\td{k}$ to 
the map $F \to \tau_{\geq -k}F,$
we get:
$$
\begin{array}{ccc}
 Rhom(\tau_{\geq -k}F, \rat_Z) & \lorw & Rhom(F, \rat_Z)  \\
\uparrow  &  & \uparrow  \\
\td{k} Rhom(\tau_{\geq -k}F, \rat_Z) & 
\stackrel{}\lorw & \td{k} Rhom( F, \rat_Z).
\end{array}
$$

\n
To prove the statement it is enough to show that
the three complexes $ Rhom(\tau_{\geq -k}F, \rat_Z),$
$\td{k} Rhom(\tau_{\geq -k}F, \rat_Z)$ and $ \td{k} Rhom( F, \rat_Z)$
have the same cohomology sheaves.
Since $F$  and $\rat_{Z}$
have locally constant cohomology sheaves, 
there are natural isomorphisms of  complexes of vector spaces 
$Rhom (F,\rat_{Z})_{y} \simeq Rhom(F_{y}, \rat_{y})
\simeq \oplus_i Hom( {\cal H}^{-i}F_{y}, \rat_{y})[-i]$.
The cohomology sheaves of the three complexes, are, therefore,
equal to $Hom( {\cal H}^{-i}F_{y}, \rat_{y})$ for $i\leq k$ and 
vanish otherwise.
\blacksquare

\medskip
The construction of the perverse truncation is done by induction on 
the strata of $Y$ starting from the shifted standard truncation on 
the 
open stratum $U_{d}.$ In the sequel we will indicate by $U_l$ the 
union
of strata of dimension bigger than or equal to $l$.
With a slight abuse of notation, we will write $U_{l+1}=U_l \coprod 
S_l,$
with $S_l$ now denoting the union of strata of dimension $l.$
Let 
$F \in Ob (D(Y))$  be  ${\frak Y}-$constructible for 
some stratification $\frak Y= \coprod S_l.$ 
All the constructions below  will lead to 
${\frak Y}-$constructible complexes. 

\medskip
We define  $\ptd{0}^{U_{d}}=\tau_{\leq -\dim{Y}}$ and
$\ptu{0}^{U_{d}}=\tau_{\geq -\dim{Y} }$.

\n
Suppose that $\ptd{0}^{U_{l+1}}: D(U_{l+1}) \lorw D(U_{l+1}) $ 
and $\ptu{0}^{U_{l+1}}: D(U_{l+1}) \lorw D(U_{l+1})$ 
have been defined. 

\smallskip
\n
We proceed to define $\ptd{0}^{U_l}$ and $\ptu{0}^{U_l}$ on 
$U_l=U_{l+1} \coprod S_l$. 
Let $i:S_l \to U_l \longleftarrow U_{l+1}: j $
be the inclusions: the exact triangles
$$
\tau'_{\leq 0}F \to F \to 
Rj_* \,^{p}\tau_{ > 0 }^{U_{l+1}}j^*F \stackrel{[1]}\lorw
\qquad
\tau''_{\leq 0}F \to F \to i_* \tau_{>-dim S } i^*F 
\stackrel{[1]}\lorw
$$
and 
$$
Rj_!\,^{p}\tau_{ < 0 }^{U_{l+1}}j^!F \lorw F \lorw \tau'_{\geq 
0}F 
\stackrel{[1]}\lorw
\qquad 
i_!\tau_{ < -\dim{S} }i^!F \lorw F \lorw \tau''_{\geq 0}F 
\stackrel{[1]}\lorw
$$
define four functors (cf. \ci{bbd},
1.1.10, 1.3.3 and 1.4.10), i.e. the four objects  
$\tau'_{\geq 0}F,$ $ \tau'_{\leq 0}F,$ $\tau''_{\geq 0}F$ 
and $ \tau''_{\leq 0}F$
which make the corresponding triangles exact, 
are determined up to  unique isomorphism.
Define
$$
\ptd{0}^{U_{l}}:= \tau''_{\leq 0} \tau'_{\leq 0}, 
\qquad \ptu{0}^{U_{l}}:= \tau''_{\geq 0} \tau'_{\geq 0}.
$$
Define:
$$
\ptd{0}:=\ptd{0}^{U_{0}}, \qquad \ptu{0}:=\ptu{0}^{U_{0}}.
$$     
We have the following compatibilities with respect to shifts.
 $$
     \ptd{m} (F [l] ) \simeq \ptd{m+l} (F) [l], \qquad
     \ptu{m} (F [l] ) \simeq \ptu{m+l} (F) [l].
     $$
These formulas hold for the ordinary truncation functors as well and 
we  symbolically 
summarize them as follows
$$
(\tau_{m} ( [l] ))[-l] \, = \, \tau_{m+l}.
$$
The perverse truncations so defined have the following properties:

\begin{itemize}
\item
By the construction above, if 
$F$ is ${\frak Y}-$cc, then so are
  $\ptd{m}F$ and $\ptu{m}F$.

\item
Let $P(Y)$ be the full subcategory of complexes $Q$ such that
 $$
\dim  \hbox{ Supp }( {\cal H}^{-i}(Q) \leq i \hbox{ for every } i \in 
\zed
$$  
and the same holds for ${\cal D}(Q),$ the Verdier dual of $Q$.
$P(Y) $ is an abelian category.
   The functor 
   $$
   \phix{0}{-}: D(Y) \lorw P(Y), \qquad \phix{0}{F}: =
   \ptd{0} \ptu{0}F \simeq \ptu{0} \ptd{0} F,
   $$
is cohomological.
   Define
   $$
   \phix{m}{F}:= \phix{0}{F[m]}.
   $$
   These functors are called {\em the perverse cohomology functors}.
   Any distinguished triangle $F \lorw G \lorw H \stackrel{[1]}\lorw$
   in $D(Y)$
   gives rise to a long exact sequence  in $P(Y)$:
   $$
    \ldots \lorw \phix{i}{F} \lorw  \phix{i}{G} \lorw
    \phix{i}{H} \lorw \phix{i+1}{F} \lorw \ldots.
   $$
If $F$ is 
${\frak Y}-$cc, then so are  $\phix{m}{F},$ $\forall m  \in \zed.$
   
\item
  Poincar\'e- Verdier Duality  induces functorial 
    isomorphisms for $F \in Ob ( D(Y)  )$
    $$\ptd{0}{\cal D}F \simeq {\cal D}\ptu{0}F, \qquad \ptu{0}{\cal 
D}F 
\simeq {\cal D}\ptd{0}F
\qquad
    {\cal D} ( \phix{j}{F})  \simeq \phix{-j}{ {\cal D} (F) }.
    $$ 
This can be seen from the construction above. In fact, 
by Lemma \ref{troncodiaz}, 
the isomorphisms hold for $U=U_{d},$ since 
$\ptd{0}^{U_d}=\tau_{\leq -\dim{Y}}$ and
$\ptu{0}^{U_d}=\tau_{\geq -\dim{Y}}$. 

\n
Suppose that 
$\ptd{0}^U{\cal D}\simeq {\cal D}\ptu{0}^U$ and
$ \ptu{0}^U{\cal D} \simeq {\cal D}\ptd{0}^U$ 
for $U=U_{l+1}$. It then follows that the same isomorphisms hold for 
$U=U_l.$
In fact, applying the functor ${\cal D}$ to the triangle defining 
$\tau'_{\leq 0}{\cal D}F$,
and the inductive hypothesis $\ptd{0}^U{\cal D}\simeq {\cal 
D}\ptu{0}^U$,
we get the triangle defining  $\tau'_{\geq 0}F$,
so that ${\cal D}\tau'_{\leq 0}{\cal D}F \simeq  \tau'_{\geq 0}F.$ 
The argument for $\tau''_{\leq 0}$ is identical.  We get
${\cal D}\tau''_{\leq 0}{\cal D}F \simeq  \tau''_{\geq 0}F.$
It follows that  
${\cal D} \tau''_{\leq 0} \tau'_{\leq 0} \simeq
\tau''_{\geq 0} {\cal D} \tau'_{\leq 0}\simeq
\tau''_{\geq 0} \tau'_{\geq 0} {\cal D}$
and the first wanted isomorphism follows.
The second is equivalent to the first one. The third
one follows formally:
${\cal D}( \phix{m}{F}) \simeq
{\cal D} \ptd{0} \ptu{0} (F[-m]) \simeq
\ptu{0}\ptd{0}  ({\cal D}(F) [m]) \simeq 
\phix{-m}{{\cal D}F}.$

\item    
     For every $F$ and $m$   one constructs,
     functorially, a distinguished triangle
     $$
      \ptd{m} F \lorw F \lorw \ptu{m+1} F \stackrel{[1]}\lorw.
   $$

\end{itemize}

 The objects of the abelian category
   $P (Y)$ are called {\em perverse sheaves.} 
 An object $F$ of $D(Y)$ is perverse if and only if
   the two natural maps $\ptd{0} F \lorw F$ and $  F \lorw  \ptu{0} F 
$
are isomorphisms.

\smallskip

\begin{ex}
\label{semismall}
{\rm Let $f:X \lorw Y$ a surjective proper map of surfaces, $X$ 
smooth. The direct image
$Rf_* \rat [2]$ is a perverse sheaf.}  
\end{ex}

\medskip

\begin{ex}
\label{3foldtrunc}
{\rm 
To give an example of how the truncation  functors can be computed 
from 
the construction given above, let us examine
the example of section \ref{lde}.
The assumptions in \ref{3foldass} are in force and we use the same 
notation.
We  show that:
$$
\ptd{0}Rf_* 
\rat_X[3]\simeq \tau_{\leq 0}Rf_* \rat_X[3],
 \qquad
 \ptd{-1}Rf_* \rat_X[3]=H_4(D)_y[1].
$$ 
Since $^{p}\!\tau^{Y-y}_{> 0}=\tau_{> -3}$ 
and $j^*Rf_*\rat_X[3]=\rat_{Y \setminus y}[3],$
we have $\tau '_{\leq 0}Rf_* \rat_X[3]= Rf_* \rat_X[3].$ 
The perverse truncation 
$\ptd{0}Rf_* \rat_X[3]=\tau ''_{\leq 0} \tau'_{\leq 0}Rf_* \rat_X[3]=
\tau ''_{\leq 0}Rf_* \rat_X[3]$
is computed by the triangle
$$
 \tau ''_{\leq 0}Rf_* \rat_X[3] \lorw Rf_* \rat_X[3] \lorw i_*\, 
 \tau_{>0}\,i^*Rf_* \rat_X[3]  
 \stackrel{+1}{\lorw}.
$$
Since 
$i^*Rf_* \rat_X[3] = \oplus_{j}  H^{3-j}(D)_y[j],$ we have 
$i_* \tau_{>0}i^*Rf_* \rat_X[3] = H^4(D)_y[-1]$, so that
$$
\ptd{0}Rf_* \rat_X[3]\, \simeq  \,
\mbox{Cone} \, \{Rf_* \rat_X[3] \to   H^4(D)_y[-1] \} \, \simeq
\, \tau_{\leq 0}Rf_* \rat_X[3].
$$
Keeping in mind the truncation rules, we have the triangle
$$
 \tau '_{\leq -1}Rf_* \rat_X[3] \lorw Rf_* \rat_X[3] \lorw Rj_* 
\,^{p}\!\tau_{>-1}^{Y-y} 
j^* \rat_Y[3] =Rj_* j^*\rat_Y[3] \stackrel{+1}{\lorw}
$$
from which we deduce that 
$$
\tau '_{\leq -1}Rf_* \rat_X[3] \,  \simeq \, 
i_!i^! Rf_* \rat_X[3] \, \simeq \, \oplus_{j} H_j(D)_y[j-3].
$$
The truncation 
$\ptd{-1}Rf_* \rat_X[3]=\tau ''_{\leq -1} \tau'_{\leq -1}Rf_* 
\rat_X[3]=
\tau ''_{\leq -1}(\oplus_{j}{ H_j(D)_y[j-3])}$
is computed by the triangle
$$
\tau ''_{\leq -1}(\oplus_{j}{ H_j(D)_y[j-3])} \lorw \oplus_{j}{ 
H_j(D)_y[j-3]} 
\to i_* \tau_{>-1}(\oplus_{j}{ H_j(D)_y[j-3] ) } \stackrel{+1}{\lorw}
$$
from which  the conclusion follows.
}
\end{ex}

\bigskip
It is remarkable  that the category of Perverse sheaves is Artinian 
and 
Noetherian,   that
its simple object can be completely characterized and  
have an important geometric 
meaning:
they are the intersection cohomology complexes. 
  
\subsection{The simple objects of $P(Y)$}
\label{simple}
Goresky and MacPherson introduced the  intersection cohomology
groups of $Y$  for an arbitrary perversity. Here we deal with the 
case 
of middle perversity. These  groups were first defined as 
the homology of a chain sub-complex of the complex of geometric
chains with twisted coefficients on $Y.$
Later, following a suggestion by Deligne, they realized these groups 
as the hypercohomology
of what they called the intersection cohomology complexes with 
twisted 
coefficients of $Y.$

\n
These complexes are the building blocks of $P(Y).$
They are special examples of perverse sheaves and every perverse sheaf
can be exhibited as a finite series of non trivial extensions
of objects of this kind supported on closed subvarieties of $Y.$

\medskip
Let $Z \subseteq Y$ be a closed subvariety,
$Z^{o} \subseteq Z_{reg} \subseteq Z $
be an inclusion of Zariski-dense open subsets and $L$
be a local system on $Z^{o}.$

Goresky-MacPherson associate with this data
the  intersection cohomology complex
$IC_{Z}(L)$ in $P(Z).$ 

\n
Up to isomorphism,  this complex
is independent of the choice
of $Z^{o}:$ if $L$ and $L'$ are local systems on
$Z^{o}$ and ${Z^{o}}'$ respectively and $L_{|Z^{o}\cap {Z^{o}}' } 
\simeq
L'_{| Z^{o}\cap {Z^{o}}' },$ then the associated intersection
cohomology complexes on $Z$ are canonically  isomorphic.

\n
The complex $IC_{Z}(L),$ when viewed as a complex on $Y,$ is perverse 
on $Y.$

\n
The intersection cohomology complex 
of $Y$ is defined to be $IC_{Y}:= IC_{Y}(\rat_{Y_{reg}}).$
If $Y$ is  smooth, or a rational homology manifold,  then 
$IC_{Y}\simeq \rat_{Y}[\dim{Y}].$

\n
If $Z$ is smooth and $L$ is a local system on $Z,$ then
$IC_{Z}(L) \simeq L[\dim{Z}].$

\begin{pr}
    \label{so}
The simple objects in $P(Y)$ are precisely the ones of the form
$IC_{Z}(L),$ $L$ simple on $Z^{o}.$
In particular, if $L$ is simple, then
$IC_{Z}(L)$  does not decompose into non-trivial
direct summands in $D(Y).$

\n
The semisimple objects of $P(Y)$ are finite direct sums
of such intersection cohomology
complexes on possibly differing subvarieties.
\end{pr}

Every perverse sheaf  $Q \in P(Y)$ is supported on a finite union
of closed subvarieties of $Y.$ Let
$Z$ be any one of them. There is a Zariski-dense 
open subset
$Z^{o}\subseteq Z_{reg},$ such that
 $Q_{|Z^{o}} \simeq L[\dim{Z}],$
where $L$ is a local system on $Z^{o}.$
The object $Q$
admits a finite  filtration where one of the quotients
is $IC_{Z}(L)$ and all the  others are all supported
on $Supp(Q) \setminus Z^{o}.$ It follows that
 $Q$ admits a finite filtration
where the quotients are intersection cohomology complexes
supported on closed subvarieties of $Y.$

An intersection cohomology complex $IC_{Z}(L)$  is characterized by 
its
not admitting subquotients supported on smaller dimensional
subspaces of $Z.$ Its eventual splitting is entirely due
to a corresponding splitting of $L.$

Let us define the intersection cohomology complexes.
Assume ${\frak Y}$ is a stratification and $L$ is a local system
on the open stratum $U_d.$ We start by defining $IC_{U_d}(L):=L[dim 
Y].$
Now suppose inductively that  $IC_{U_{l+1}}(L)$ has been defined on 
$U_{l+1}$
and we define it on $U_l$ by 
$$
IC_{U_l}(L):=\tau_{\leq -l-1}Rj_*IC_{U_{l+1}}(L).
$$

\bigskip
Let us give formulae for $IC_{Y}(L)$ when $Y$ and $L$ have isolated 
singularities. It suffices to work in the Euclidean topology.

\n
Let $(Y,p)$  be a germ of an isolated singularity,
$j:  U: = Y \setminus p  \to Y$ be the open embedding 
and $L$ be a local system on $U.$
We have
\begin{equation}
    \label{icis}
IC_{Y}(L) = \td{-1} (Rj_{*}L[\dim{Y}]).
\end{equation}
If $\dim{Y} =1,$ then
$IC_{Y}(L) =  j_{*}L[1].$
The stalk at $p$ are the invariants of $L$ i.e.
$H^{0}(U, L).$

\n
In general, when  $\dim{Y} \geq 2,$ then
$IC_{Y}(L)$ is a complex, not a sheaf.
If $L$ is simple, then $IC_{Y}(L)$ is simple
and does not split non-trivially in $D(Y).$
The 
cohomology sheaves ${\cal H}^{j}(IC_{Y}(L))$
are non trivial only for $j \in [-\dim{Y}, -1]$ and we have
\begin{equation}
    \label{csisl}
{\cal H}^{-\dim{Y}}(IC_{Y}(L)) = j_{*}L, \qquad
{\cal H}^{-\dim{Y}+l}(IC_{Y}(L)) = H^{l}(U,L)_{p}, \;  1 \leq  l \leq
\dim{Y} -1,
\end{equation}
where $V_{p}$ denotes a skyscraper sheaf at $p\in Y$
with stalk $V.$

\medskip

In order to familiarize ourselves  with these complexes, we  compute 
two important examples:

\begin{ex} 
\label{3folddoublepoint}
{\rm We consider a threefold $Y$ with an ordinary double point $y$
and with associated link  ${\cal L}.$ 
Let $j:Y\setminus y \to Y$ be the open embedding, so that 
(\ref{icis}) gives 
$$
IC_Y= \tau_{\leq -1}Rj_*\rat[3].
$$ 
The cohomology sheaves at $y$ are
$$
{\cal H}^{k}(IC_{Y}) =  H^{k+3}({\cal L})_{p}, \hbox{ for } k 
\leq -1, \qquad {\cal H}^{k}(IC_{Y})= 0 
\hbox{ otherwise. } 
$$
The singularity is analytically equivalent to a cone over a smooth 
quadric 
in projective space,
hence its link is homeomorphic to the $S^1-$bundle over $S^2 \times 
S^2$ 
with Chern class $(1,1).$
The long exact sequence for this $S^1-$fibration gives 
$$
H^k({\cal L})= \rat \hbox{ for }  k=0,2,3,5 \qquad H^k({\cal L})= 0 
\hbox{ otherwise, }  
$$
which in turn implies that 
$$
{\cal H}^{k}(IC_{Y}) =  \rat  \hbox{ for } k= -3, -1,  \qquad
{\cal H}^{k}(IC_{Y}) = 0 \hbox{ otherwise. } 
$$
We have a triangle in $D(Y),$ (not of perverse sheaves)
$$
\rat_Y[3] \to IC_Y \to \rat_y[1] \stackrel{+1}{\to}
$$
The fact that 
${\cal H}^{-1}(IC_{Y}) =  \rat_y $ should be compared with the 
existence 
of a small resolution 
with fiber a projective line over the singular point, and the 
statement 
of the Decomposition Theorem \ref{dtbbd}.}
\end{ex} 
\begin{ex} 
\label{cksrecipe}
{\rm 
Let $Y=\comp^2,$ and $L$ be a local system on $\comp^2 \setminus ( 
x_1x_2=0)$ 
defined by the two monodromies $T_1$ and $T_2$ acting 
on the vector space $V=L_{p},$ the stalk of $L$ at $p=(1,1) \in 
\comp^{2}.$ 
We first determine the intersection cohomology complex
over $\comp^2 \setminus \{0\}.$ 
Denoting by  $j:\comp^2 \setminus \{ x_1x_2=0\} \to \comp^2 \setminus 
\{0\}$ 
the natural map, we have 
 $IC_{\comp^2 \setminus \{0\} }(L)=\tau_{\leq -2}Rj_*L[2]=(j_*L) [2].$
Denoting by 
$j':\comp^2 \setminus \{0\} \to \comp^2$  the natural map, we have  
$$
IC_{\comp^2}(L)\, = \, \tau_{\leq -1}Rj'_* IC_{\comp ^2 \setminus 
\{0\} (L) }  
\, =\, \tau_{ \leq -1}Rj'_*(j_*L [2]).
$$
In order to determine  the cohomology sheaves of $IC_{\comp^2}(L),$ 
we compute
$H^i(\comp^2 \setminus \{0\}, j_*L)$ for $i=0,1.$
More precisely, we should determine these groups for a fundamental 
system 
of neighborhoods of the origin; however  the cohomology groups are in 
fact constant.
Set $N_1=T_1-Id,$ $N_2=T_2-Id.$

\smallskip
\n
We have  
$H^0(\comp^2 \setminus \{0\}, j_*L)=H^0(\comp^2 \setminus \{ 
x_1x_2=0\},L)
= \ke{ N_1} \cap \ke{ N_2}=V^{\pi_1} .$ 
Since $j_*L= \tau_{\leq 0}Rj_*L,$ and fundamental deleted  
neighborhoods 
around
the axes are homotopic to circles, so that
$$
{\cal H}^{i}(Rj_*L) =0 \hbox{ for } i \geq 2,  
$$ 
we have the following exact triangle in $D(\comp^2 \setminus \{0\})$
 $$
j_*L =\tau_{\leq 0}Rj_*L \lorw Rj_*L \lorw   
{\cal H}^{1}(Rj_*L)[-1] \stackrel {+1} {\lorw}. 
$$
The sheaf
${\cal H}^{1}(Rj_*L)$ is the local system on  
$( x_1x_2=0) \setminus \{(0,0) \}=D_1 \coprod D_2$,
with fiber $\coker{ N_1}$ and monodromy $T_2$ on $D_1$, and 
fiber $\coker{ N_2}$ and monodromy $T_1$ on $D_2$. 
Since $\comp^2 \setminus \{ x_1x_2=0\}$ retracts to a torus $T^2$, 
the cohomology of $L$ is isomorphic to the group cohomology of $\zed 
^2$
with values in $V$ as a $\zed ^2$-module via the monodromies 
$T_1,T_2,$
which can be computed by the Koszul complex (see for instance
\ci{weibel}                                                              )

$$
0 \lorw V \stackrel{\phi}{\lorw} 
 V \oplus V \stackrel{ \psi}{\lorw}V {\lorw} 0, 
$$
with 
$$
\phi(v)= (N_1(v),N_2 (v)) \qquad \psi(v_1,v_2)=N_2(v_1) -N_1(v_2).
$$
The long exact sequence associated to the exact triangle above gives
$$
{\cal H}^{-1}(IC_Y(L))_0 \simeq
H^1(\comp^2 \setminus \{0\}, j_*L) \,= \,
\frac{(N_1(v_1),N_2(v_2))\hbox{ such that }N_1N_2(v_1-v_2)=0 }
{ (N_1( v),N_2( v))}.
$$
More generally, a similar recipe holds for the cohomology sheaves 
of the intersection cohomology complex of a local system 
defined on the complement of a normal crossing, see \ci{cks}.}
\end{ex}

\subsection{Decomposability,  $E_{2}-$degenerations
and filtrations}
\label{de2}
\begin{defi}
\label{hdec}
{\rm
Let ${\cal H}= {\cal H}^{0}$ be the sheaf cohomology functor.
We say that $F$ in $D(Y)$ is {\em ${\cal H}-$decomposable} if
$$
F \simeq_{D(Y)} \bigoplus_{i}{ {\cal H}^{i}(F)[-i] }
$$
We say that $F$ in $D(Y)$ is {\em $^{p}{\cal H}-$decomposable} if
$$
F \simeq_{D(Y)} \dsdix{i}{F}.
$$
}
\end{defi}

If $F$ is ${\cal H}-$decomposable, then the spectral sequence
$$
{H}^{p}( Y, {\cal H}^{q}(F) ) \Longrightarrow {\Bbb H}^{p+q}(Y,F)
$$
is $E_{2}-$degenerate.
This spectral sequence is the Leray Spectral Sequence when
$F= Rf_{*}(G).$ In this case the corresponding filtration
is called the Leray filtration.

\n
The analogous statements holds for $^{p}\!{\cal H}-$decomposability.
The corresponding spectral sequence is called the Perverse
Leray Spectral Sequence:
$$
{\Bbb H}^{p}(Y,  ^{p}{\cal H}^{q}(Rf_{*}G) ) \Longrightarrow {\Bbb 
H}^{p+q}(Y,G)
$$
and the corresponding filtration is called the perverse filtration.
\begin{defi}
\label{pf}
{\rm
 Let $f: X \to Y$ be a map, $n = \dim{X}.$
 The perverse filtration  $H^{n+j}_{\leq b}(X) \subseteq H^{n+j}(X),$ 
$ 
 b,j \in \zed$
 is defined to be the perverse filtration
 on ${\Bbb H}^{j}(Y, Rf_*\rat_{X}[n])$ 
 }
 \end{defi}

It coincides, up to a shift, with the Leray filtration, when $f$ is 
smooth.

\medskip
If these decomposing isomorphisms exist, they are seldom unique.
We now give the statement (not in the most general form)
of one of the more general criteria for decomposability,
see \ci{dess} and \ci{shockwave}:

\begin{tm}
\label{delilefsdege}
{\bf \rm (Deligne degeneration criterion.)} 
Let $K$ be an object of $D^b(Y)$, and let 
$\eta \in H^2(X).$ Suppose that $\eta^l :  
\phix{-l}{K} \to \phix{l}{K}$ 
is an isomorphism for all $l.$
Then $K$ is $^{p}\!{\cal H}-$decomposable.
The same statement holds if we consider the functor ${\cal H}$.
\end{tm}

\begin{ex}
\label{3foldancor}
{\rm By the computation done in \ref{3foldtrunc},
we have the following description of the perverse filtration 
for the resolution of a threefold:
$$
H^{i}_{\leq -2}(X)=\{0\}, \qquad H^2_{\leq -1}(X)=\im \{H_4(D) \to 
H^2(X)\}, \qquad 
H^i_{\leq -1}(X)=0 \hbox{ otherwise, }
$$
$$
H^4_{\leq 0}(X)=\ke \{H^4(X) \to H^4(D)\}, \qquad 
H^i_{\leq 0}(X)=H^i(X) \hbox{ otherwise }, 
$$
$$
 H^i(X)_{\leq 1}=H^i(X) \hbox { for all } i.
$$
The condition \ref{delilefsdege},  that
$$\eta:  
\phix{-1}{Rf_* \rat_X[3]}=
H_4(D)_y \lorw H^4(D)_y=\phix{1}{Rf_* \rat_X[3]}
$$
be an isomorphism,  is just        
the non degeneracy of the intersection form 
$\eta \cap [D_i]\cdot[D_j].$
Note that in this case, the explicit description makes it clear that
the perverse filtration is given by Hodge substructures.}
\end{ex}

\subsection{The Decomposition Theorem of
Beilinson, Bernstein, Deligne and Gabber}
\label{sdtbbd}
We can now state the  generalization of 
Deligne's Theorem \ref{del}
to the case of arbitrary proper maps.
Recall that if $X$ is smooth, then $IC_{Y} = \rat_{X}[n].$
\begin{tm}
\label{dtbbd}
Let $f:X \to Y$ be a proper map of algebraic
varieties. Then
\begin{equation}
    \label{fd}
    \mbox{the complex $Rf_{*}IC_{X} \simeq
    \dsdix{i}{Rf_{*}IC_{X}}$ is $^p\!{\cal H}-$decomposable }
    \end{equation}
The complexes 
$\phix{j}{Rf_{*}IC_{X}}$ are semisimple i.e.
there is a canonical isomorphism
\begin{equation}
    \label{semisi}
 \phix{j}{Rf_{*}IC_{X}} \,\simeq_{P(Y)}  \,
 \oplus{IC_{Z_{a}}(L_{a} )}
\end{equation}
for some finite collection, depending on $j,$
of semisimple local systems $L_{a}$ on smooth 
\underline{distinct} varieties $Z_{a}^{o} \subseteq Z \subseteq Y.$

\smallskip
\n
Let $\eta$ be an $f-$ample line bundle on $X.$
 Then
\begin{equation}
    \label{rhl}
\eta^{r} \, :  \,  \phix{-r }{  Rf_{*}IC_{X} } \,  \simeq \,
\phix{r }{  Rf_{*}IC_{X} }.
    \end{equation}
\end{tm}

\medskip
The Verdier Duality functor  is an autoequivalence ${\cal D}: D(Y) 
\to D(Y)$ which preserves $P(Y)$ and for which
one has 
$$
{\cal D} \circ \, ^{p}\!{\cal H}^{-j} \, \simeq \, 
^{p}\!{\cal H}^{j} \circ {\cal D}.
$$
This fact
implies that
the summands appearing in the semisimplicity statement 
for $j$ are pairwise isomorphic to the ones appearing
for $-j$ and that the local systems $L$ are self-dual.

\medskip
Theorem \ref{dtbbd} is the deepest known fact concerning the homology
of algebraic maps.

\noindent
The original proof uses algebraic geometry in positive characteristic.
in an essential way.

\n
M. Saito has given a transcendental proof of a more general statement
concerning his mixed Hodge modules in the series of papers 
\ci{samhp}, \ci{samhm}, \ci{samhp}.

\n
  We give a proof for the 
push-forward of intersection cohomology (with constant coefficients)
first in  the case of semismall maps (cf. \ci{demigsemi}) and
 then for arbitrary maps in  \ci{decmightam}. Though at present 
 our  methods do not  afford results concerning the push-forward with 
 more general coefficients,  they give new and precise results on the 
perverse filtration and on the refined intersection forms.

\n
C. Sabbah \ci{sa} has recently proved a decomposition theorem for 
push-forwards of semisimple local systems.

\begin{rmk}
\label{infatt}
{\em It is now evident that the computations in \ref{fsts} and 
\ref{fstc}
establish the Decomposition Theorem for maps from a smooth surface.
In the case of the proper birational map $f:X \to Y$ of \ref{fsts}, 
in fact, 
the complex $Rf_* \rat_X[2]$ is perverse, as observed in 
\ref{semismall}, and \ref{tmrs}
states that it splits into $IC_Y$ and $R^2f_* \rat_X[0].$
In the case of the family of curves treated in \ref{fstc} we have that
$j_*T^0[2]=\phix{-1}{  Rf_{*}\rat_{X}[2] }[1],$ and 
$j_*T^2[0]=\phix{1}{  Rf_{*}\rat_{X}[2] }[-1],$ and we 
showed in \ref{rhl1} that
$\eta: j_*T^0 \to j_*T^2$ is an isomorphism.
The perverse sheaf $P$ splits, see \ref{semisimp1}, in 
$j_*T^1[1]=IC_Y(T^1),$ 
and $V$, concentrated on points. }
\end{rmk}

\begin{rmk}
\label{infatt3fold}
{\em For the case of  the 
resolution of a threefold with isolated singularities,
whose Hodge theory has been  
treated in \ref{lde}, we have, 
as seen in \ref{3foldtrunc}, 
$$\phix{-1}{  Rf_{*}\rat_{X}[3] }\simeq H_4(D)_y, \qquad 
\phix{1}{  Rf_{*}\rat_{X}[3] }\simeq H^4(D)_y \simeq \eta \wedge 
H_4(D)_y,$$
and we have the splitting
$$
\phix{0}{  Rf_{*}\rat_{X}[3] }\simeq IC_Y \oplus H_3(D)_y.
$$
Similarly, for the 4-fold with isolated singularities, see \ref{lde4},
$$\phix{-2}{  Rf_{*}\rat_{X}[4] }\simeq H_6(D)_y, \qquad 
\phix{-1}{  Rf_{*}\rat_{X}[4] }\simeq H_5(D)_y ,$$
$$\phix{2}{  Rf_{*}\rat_{X}[4] }\simeq H^6(D)_y\simeq \eta^2 \wedge 
H_6(D)_y, \qquad 
\phix{1}{  Rf_{*}\rat_{X}[4] }\simeq H^5(D)_y \simeq \eta \wedge 
H_5(D)_y,$$
and we have the splitting 
$$
\phix{0}{  Rf_{*}\rat_{X}[4] }\simeq IC_Y \oplus H_4(D)_y.
$$}
\end{rmk}

\subsection{Results on intersection forms }
\label{inr}

In this section we list some of the results of
\ci{decmightam} which are related to the theme of this paper.
For simplicity, we state them in the special case when $f: X \to Y$
is a map of projective 
varieties, $X$ smooth. Let $\eta$ and $A$
be ample line bundles on $X$ and $Y$ respectively, $L:= f^{*}A.$

\begin{tm}
\label{uf}
For $l\geq 0$ and $b\in \zed,$ the subspaces given by the perverse 
filtration (cf. \ref{de2})
$$
H^{l}_{\leq b}(X) \, \subseteq\,  H^{l}(X)
$$
are  pure Hodge sub-structures. The quotient spaces
$$
H^l_b(X)\, = \, H^{l}_{\leq b}(X) /H^{l}_{\leq b-1}(X)
$$
inherit a  pure Hodge structure 
of  weight $l.$
\end{tm}

The cup product with $\eta$ verifies 
$\eta \, H^l_{\leq a}(X) \subseteq H^{l+2}_{\leq a+2}(X)$ 
and induces maps,  still denoted 
$\eta: H^l_{a}(X) \to H^{l+2}_{a+2}(X)$.
The cup product with $L$ is  compatible with the Decomposition Theorem
\ref{dtbbd}
and induces maps $L: H^l_{a}(X) \to H^{l+2}_{a}(X).$

\n
These maps satisfy  graded  Hard Lefschetz Theorems (cf. 
\ci{decmightam}, Theorem 2.1.4).

\n
Define 
$P^{-j}_{-i}:= \ke{\, \eta^{i+1}} \cap \ke{\, L^{j+1}} 
\subseteq H^{n-i-j}_{-i}(X),$ 
$i,\,j \geq 0$ and   $P^{-j}_{-i} :=0$
otherwise.
In the same way in which the classical Hard Lefschetz implies
the Primitive Lefschetz Decomposition for the cohomology of $X,$
the graded Hard Lefschetz Theorems  imply the 
double direct sum decomposition of

\begin{tm}
\label{etaldecompo}
Let $i, \, j \in \zed.$ 
There is  a Lefschetz-type direct sum decomposition into 
pure Hodge sub-structures  of weight $(n-i-j),$
called the  $(\eta,L)-$decomposition:
$$
H^{n-i-j}_{-i}(X) = \bigoplus_{l,\,m\, \in \zed }{
\eta^{-i+l} \, L^{-j+m} \, P^{j-2m}_{i-2l}.}
$$
\end{tm} 

One can define bilinear forms 
$S^{\eta L}_{ij}$ on $H^{n-i-j}_{-j}(X)$ by 
modifying the Poincar\'e pairing
$$
S^{\eta L}_{ij} ([\alpha], [\beta]   ) \, := \,
\int_X{ \eta^i \wedge  L^j \wedge  \alpha \wedge \beta}
$$
and descending it to the graded groups.
These forms are  non degenerate. In fact their signature
can be determined in the following generalization of the Hodge
Riemann relations.
 
\begin{tm}
\label{tmboh}
The \ref{etaldecompo} is orthogonal with respect
to $S^{\eta L}_{ij}.$ The forms $S^{\eta L}_{ij}$
induce polarizations 
of each $(\eta,L)-$direct summand.
\end{tm}

The homology groups $H^{BM}_{*}(f^{-1}(y))=H_{*}(f^{-1}(y)),$
$y \in Y,$ 
are filtered by virtue of the decomposition theorem (one may call this
the perverse filtration).
The natural cycle class map
$cl: H^{BM}_{n-*}(f^{-1}(y)) \to H^{n+*}(X)$ is filtered strict.

\medskip
The following generalizes the Grauert Contractibility Criterion. 

\begin{tm}
\label{nhrbr} 
Let $b \in \zed$, $y \in Y$.
 The  natural class  maps
$$
cl_{b}: H^{BM}_{n-b,b}(f^{-1}(y)) \lorw H^{n+b}_b(X)
$$
is injective and identifies
$H^{BM}_{n-b,b}(f^{-1}(y))\subseteq \ke{\, L}\subseteq
H^{n+b}_b(X)$ with a pure  Hodge substructure, compatibly with
the $(\eta,L)-$decomposition. Each $(\eta,L)-$direct summand
of $H^{BM}_{n-b,b}(f^{-1}(y))$ is polarized up to sign
by $S^{\eta L}_{-b,0}.$

\n
In particular, the restriction of $S^{\eta L}_{-b,0}$ to
$H^{BM}_{n-b,b}( f^{-1} (y ))$ 
is non degenerate.
\end{tm}

By intersecting in $X$ cycles supported on $f^{-1}(y),$
we get the refined intersection form (see section  \ref{psif})
$H^{BM}_{n-*}(f^{-1}(y)) \to H^{n+*}(f^{-1}(y))$ 
which is filtered strict as well.

\begin{tm}
\label{rcffv}
({\bf The Refined Intersection Form Theorem})
Let $b \in \zed$, $y \in Y$.
The graded refined intersection form
$$
H^{BM}_{n-b,a}(f^{-1}(y)) \lorw H^{n+b}_{a}(f^{-1}(y))
$$
is zero if $a\neq b$ and it is an isomorphism if $a=b.$
 \end{tm}

We have seen in earlier sections how these results can be made 
explicit in the case of 
surfaces, threefolds and fourfolds.
For more applications in any dimension see \ci{decmightam}.

\medskip
In fact, the method of proof of the results stated
in this section  is inspired by the low
dimensional examples
of surfaces, threefolds and fourfolds. 

\subsection{The decomposition mechanism}
\label{decmec}
It is quite hard to describe what kind of geometric phenomena
are expressed by the Decomposition Theorem.
The complex $Rf_* \rat_X$ essentially
describes $H^*(f^{-1}(U))$ for any neighborhood $U$ of a point $y.$
We gain some geometric insight if we represent, via Poincar\'e
duality, the cohomology classes with  Borel-Moore cycles in $f^{-1}(U).$
Let $S$ be the stratum containing $y.$ By the exact sequence
$$
H_*^{BM}(f^{-1}(U \cap S)) \stackrel{i_*}{\to}
H_*^{BM}(f^{-1}(U)) \stackrel{j^*}{\to}
H_*^{BM}(f^{-1}(U \setminus S)), 
$$
the Borel Moore cycles in $f^{-1}(U)$ are of two kind:
those in $\im i_*$,  which are homologous to cycles supported 
on the inverse image of the stratum $S,$ and those whose restriction to 
$f^{-1}(U \setminus S)$ is not trivial.
Neither $i_*$  is injective nor $j^*$ is surjective:
there are non trivial cycles in $f^{-1}(U \cap S)$
which become homologous to zero in $f^{-1}(U )$, and there are
cycles in $f^{-1}(U \setminus S)$ which  cannot be closed 
to cycles in $f^{-1}(U).$

The Decomposition Theorem gives strong information on both types.
The first deep aspect of  the Theorem is that the subspace 
$\im i_*,$ has a uniform behavior for all projective maps,
related to the non degeneracy of the intersection forms. 
For instance, we already noticed in \ref{link}, 
how the Grauert Theorem \ref{tmgra} implies that the classes
of disks transverse to exceptional curves are homologous to 
linear combinations of the classes of these curves.
Such non degeneracy results, see \ref{nhrbr}, \ref{rcffv},
are peculiar to algebraic maps and stem from "weight'' considerations,
either in characteristic $0$ (Hodge Theory). or in positive characteristic
(weights of Frobenius, cf. \ci{bbd}).

The Decomposition Theorem, though, contains other deep information.
Since   $Rf_* \rat_X$ splits as a direct sum of terms associated 
with the strata, we have a  splitting  map, which can be made canonical after
an ample line bundle $\eta$ on $X$ has been chosen, from the subspace 
$\im{ j^* } \subseteq H_*^{BM}(f^{-1}(U \setminus S)$ to 
$H_*^{BM}(f^{-1}(U ):$ i.e. the following is split exact
$$
0 \lorw \im{i_*} \lorw H_*^{BM}(f^{-1}(U ) \lorw \im{j^*} \lorw 0.
$$
 The image of this map defines a subspace
of $H_*^{BM}(f^{-1}(U)$ which is complementary to $\im i_*$
and consists of classes which are closures of some Borel Moore cycles
in $f^{-1}(U \setminus S).$
The deep fact here, is that these cycles are governed by the 
intersection cohomology complex construction on $Y;$ each stratum having $S$ 
in its closure contributes to $H_*^{BM}(f^{-1}(U)$ via the  intersection 
cohomology of a local system on the stratum.

\section{Grothendieck motive decomposition for maps of threefolds}
\label{gmdmt}
We assume we are in the situation \ref{3foldass}.
Again, this is for ease of exposition only. See Remark
\ref{rmkfi}.

\medskip
We have already shown that
$$
H^{2}_{-1}(X) \, = \, \im{  H_{4}(D) \to H^{2}(X)    }  \, =: \, 
\im{i_{*}},
$$
$$
H^4_{0}(X)=\ke{ \{H^4(X) \to H^4(D)\}} \, =: \, \ke{ i^*}. 
$$
The choice of an ample line bundle $\eta$ allows to split 
the perverse filtration: 
$$H^2(X)=\im i_* \oplus (c_1( \eta) \wedge \im i_*)^{\perp}.$$
$$H^4(X)=\ke i^* \oplus ( \im i_*)^{\perp}.$$
We have, canonically, that 
$$H^3(X)=\im{ H_3(D)} \oplus (\im {H_3(D))}^{\perp}.$$
so that 
$$
I\!H^i(Y)\,=\,H^i(X) \hbox{ for } i=0,\,1, \,5,\,6, 
\qquad I\!H^2(Y)\, = \, (c_1( \eta) \wedge \im H_4(D),
)^{\perp}
$$
$$
I\!H^3(Y)\, = \, (\im{H_3(D)})^{\perp}, \qquad I\!H^4(Y)\, =\, (\im 
H_4(D))^{\perp}; 
$$
here we are using the convention for intersection cohomology
compatible with singular cohomology: 
$I\!H^i(Y):= {\Bbb H}^{i-n}(Y, IC_Y)).$

\bigskip
We want to realize these splittings by algebraic cycles 
on $X \times X$, in order to find a Grothendieck motive
for the intersection cohomology of $Y.$
These cycles will be supported on $D \times D.$

\medskip
We start with the following simple lemma. 
\begin{lm}
\label{support}
Let $X$ be a projective $n-$fold, 
and $Y \subseteq X$ be  a subvariety.
Let $W \subseteq 
\im{ \{ H_s(Y) \to H^{2n-s}(X) \} } \subseteq H^{2n-s}(X)$
be a vector subspace
on which the restriction  of the Poincar\`e pairing remains
non degenerate, i.e.
$H(X)=W \oplus W^{\perp}.$ 
Then the projection 
$P_W \in End(H(X))\simeq H(X \times X)$ on $W$ 
relative to the above splitting can be represented by a cycle 
supported on 
$Y \times Y$.
\end{lm}
{\em Proof. }
Let $\{e_1\}$ be a basis for $H(X)$
such that $ e_1, \cdots, e_k \in W$ and 
$ e_{k+1}, \cdots, e_N \in W^{\perp}.$
For $i=1, \cdots, k,$ we can represent $e_i$ by a cycle $\gamma_i$
contained in $Y$. In force of the hypothesis, the dual basis 
$\{ {e_i}\check{} \}$ is of the form 
$$
e_i \check{}= \sum_{j=1}^k a_{ij}e_j \qquad \hbox{ for }1 \leq i \leq 
k
\qquad
e_i \check{}= \sum_{j=k+1}^N a_{ij}e_j \qquad \hbox{ for }k+1 \leq i 
\leq N.
$$
In particular 
$e_1 \check{}, \cdots,e_k \check{} $ are represented by the cycles 
$
\gamma_i \check{}= \sum_{j=1}^k a_{ij}\gamma_j 
$ supported on $Y$. The projector 
$P_W=  \sum_{i=1}^k e_i \otimes e_i \check{}$ is thus represented by 
the cycle
$  \sum_{i,j=1}^k a_{ij}\gamma_i \times \gamma_j$, which is supported 
on
$Y$.
\blacksquare

\medskip
The following is a standard but very useful application of 
``strictness'' in Hodge Theory

\begin{lm}
\label{pesi}
Let $Y \subseteq X$ be a codimension $d$ 
subvariety of an $n-fold$
and let $\pi:\tilde{Y} \to Y$ a resolution of singularities.
Suppose $\beta \in \im \{H_{2k}(Y) \to H^{2(n-k)}(X)\} 
\cap H^{n-k,n-k}(X).$ Then there is 
$\tilde{\beta}\in  H^{n-k-d,n-k-d}(\tilde{Y})$ such that 
$(i\circ \pi)_*(\tilde{\beta})= \beta.$ 
\end{lm}
{\em Proof.} We consider the weights of the homology groups as given 
by 
their being dual of the cohomology groups.
Thus $H_{2k}(Y)$ has weights $\geq -2k.$
The map $H_{2k}(Y) \to H^{2(n-k)}(X)$ is of type $(n,n)$. Since the 
Hodge 
structure on $ H^{2(n-k)}(X)$ is pure, 
the strictness of maps of Hodge structures implies that
$$
\im{ \{ H_{2k}(Y) \to H^{2(n-k)}(X) \} } \,
= \, \im{ \{ W_{-2k}H_{2k}(Y) \to H^{2(n-k)}(X) \} }.
$$
It follows that $\beta =i_*\beta '$ for some
$\beta ' \in  W_{-2k}H_{2k}(Y).$ On the other hand this group
coincides with 
$ \im{ \{ \pi_* : H_{2k}(\tilde{Y}) \to H_{2k}(Y) \}}$
 for any resolution, whence the statement.
\blacksquare

\begin{tm}
\label{cicli}
Let $f:X \to Y$, $D$ as before.
Then there exist algebraic 3-dimensional cycles $Z_{-1}, Z_{0}, Z_1,$
supported on $D \times D$ such that:

\smallskip
\n
$Z_1$ defines the projection of $H(X)$ onto 
$H^4_1(X)= c_1( \eta) \wedge \im \{H_4(D) \to H^2(X)\} \subseteq 
H^4(X);$

\n
$Z_{-1}$ defines the projection of $H(X)$ onto 
$H^2_{-1}(X)=   \im \{H_4(D) \to H^2(X)\} \subseteq H^2(X);$

\n
$Z_0$ defines the projection of $H(X)$ on 
$ \im \{H_3(D) \to H^3(X)\} \subseteq H^3(X).$
\end{tm}
{\em Proof.}
Let $\Lambda$ be the inverse of the negative-definite intersection 
matrix 
$I_{ij}= \int_X c_1(\eta) \wedge[D_i]\wedge [D_j].$
We denote by $\eta \cap D_i$ the curve obtained intersecting the 
divisor $D_i$ 
with a general section of $\eta$. 
Set:
$$
Z_{-1}= \sum \lambda_{ij}[(\eta \cap D_i) \times D_j] 
\qquad Z_{1}= \sum \lambda_{ij}[D_i \times (\eta \cap D_j)].
$$
It is immediate to verify that $Z_{-1} $ and $Z_1$ define the 
sought-for projectors.

\n
The construction of  $Z_0$ is not so direct: 
 Since, by \ref{polarizeh3}, the Poincar\'e paring is non degenerate 
on 
$ \im \{H_3(D) \to H^3(X)\}$, by \ref{support} we can represent 
the projection on $H_3(D)$ by a cycle supported on $D \times D$.
Furthermore, the projection is a map of Hodge structures, hence its 
representative cycle 
$P_3 \in H^6(X \times X)$
has type 
$(3,3).$ By \ref{pesi} we have $P_3=i_* \pi_* \beta$ for some 
$\beta \in H^{1,1}(\widetilde{D \times D}),$ where $\widetilde{D 
\times D}$ 
is any resolution of 
${D \times D}.$
By the Lefschetz Theorem on $(1,1),$ classes there is an algebraic 
cycle $\tilde{Z}$ such
that $\beta =[\tilde{Z}]$. It is clear that $Z_0=i_*  \pi_* Z$ does 
the job.
\blacksquare

\medskip
The following follows immediately.
\begin{cor}
    \label{gm}
The Grothendieck motive, 
 $(X, \Delta_{X} -Z_0 -Z_1 -Z_{-1})$
is a  Betti realization of the Intersection cohomology of $Y.$
\end{cor}

\medskip
We can be more specific: the projector  $\Delta_{X} -Z_0 -Z_1 -Z_{-1}$
is supported on the fiber product $X \times_Y X$, 
therefore defines a relative motive over $Y$
in the sense of \ci{ch} (see also \ci{decmigmot}).
By \ci{ch}, Lemma 2.23, the isomorphism of algebras  
$$
End(Rf_* \rat_X [3])= H^{BM}_6(X \times_Y X)
$$
ensures that
the Betti realization of this relative motive is the projector 
associated with  the  splitting  for $Rf_* \rat_X [3]$ we have used 
in this section.

\begin{rmk}
{\rm From the construction of the cycles it is evident that $Z_{-1}$
and $Z_1$ define in fact Chow motives, not only Grothendieck motives.
Under some hypothesis it is possible to construct a Chow projector
for $Z_0$ as well. For instance, if $D$ is smooth irreducible, and its
conormal bundle ${\cal I}_D/{\cal I}_D^2$ is ample. In this case,
let $Z_0$ the
cycle in $D \times D$ representing the Hodge $\Lambda$ operator with
respect to the polarization given by the conormal bundle. It is
immediate to verify that $Z_0$ defines the  Chow motive we need. In
general some knowledge of the nature of the resolution may allow one
to find a Chow motive whose Betti realization is intersection 
cohomology.} 
\end{rmk}

\begin{rmk}
    \label{rmkfi}
{\rm It is not difficult to modify the proofs to produce a
Grothendieck motive for the  intersection cohomology of $Y$ for an 
arbitrary
 three-dimensional variety $Y,$ (e.g. with non-isolated 
singularities).
If, for example, 
some divisor $D'$ is blown down to a curve $C,$ then  one needs to 
construct a
further projector, represented by a cycle which is a linear
combination of the components  of $D' \times_C D'.$  This projector
 splits off the
contribution of $D'$ to the cohomology of $X$. We leave this task to
the reader.} 
\end{rmk}

\begin{rmk}
\label{chissa}
{ \rm If $Y$ is a fourfold with isolated singularities, then the 
computations in 
\ref{lde4} express its intersection cohomology as a  Hodge 
substructure
of the cohomology of a resolution $X.$ The method developed in this 
section does not apply in general since we do not know whether the 
classes of the projectors, which are pushforward of classes of type 
$(p,p)$ on a resolution of the product of the exceptional divisor 
with itself, are represented by algebraic cycles.
On the other hand, 
this can be achieved in the presence of supplementary 
information on the singularities of $Y$ or on the exceptional 
divisor. For example:    if the singularities are locally isomorphic to toric 
singularities. This allows to define a motive for the intersection 
cohomology in several interesting cases.}
\end{rmk}

Authors' addresses:

\smallskip
\n
Mark Andrea A. de Cataldo,
Department of Mathematics,
SUNY at Stony Brook,
Stony Brook,  NY 11794, USA. \quad 
e-mail: {\em mde@math.sunysb.edu}

\smallskip
\n
Luca Migliorini,
Dipartimento di Matematica, Universit\`a di Bologna,
Piazza di Porta S. Donato 5,
40126 Bologna,  ITALY. \quad
e-mail: {\em migliori@dm.unibo.it}

\end{document}